\documentclass[]{scrartcl}
\usepackage[utf8]{inputenc}
\usepackage[letterpaper,bindingoffset=0cm,inner=2.5cm,outer=2.5cm,top=2.5cm,bottom=2.5cm]{geometry}
\usepackage{helvet}
\usepackage[T1]{fontenc}
\usepackage{graphicx}
\usepackage{amsmath,amsthm,amstext,amssymb,bm}
\usepackage{mathtools}
\usepackage{xcolor,color}
\usepackage{enumerate}
\usepackage{pgf}
\usepackage[toc,page]{appendix}
\usepackage[colorlinks=true]{hyperref}
\definecolor{mylinkcolor}{RGB}{0,0,130}
\hypersetup{colorlinks,allcolors=mylinkcolor,citecolor=mylinkcolor}
\usepackage{url}
\usepackage[capitalize,nameinlink]{cleveref}
\crefformat{equation}{(#2#1#3)}
\usepackage{booktabs}
\usepackage{tabularx}
\usepackage{enumerate}
\usepackage[shortlabels]{enumitem}
\usepackage{multirow}
\usepackage{subfigure}
\usepackage[font=small,labelfont=bf]{caption}

\usepackage{import}
\usepackage{float}
\usepackage{tikz}
\usepackage{pgfplots}
\usepackage{blkarray,bigstrut}
\definecolor{DarkGreen}{rgb}{0, 0.392, 0}
\newcommand{\mLabel}[1]{\mbox{$\scriptstyle{#1}$}}
\newcommand\topstrut[1][1.2ex]{\setlength\bigstrutjot{#1}{\bigstrut[t]}}
\newcommand\botstrut[1][0.6ex]{\setlength\bigstrutjot{#1}{\bigstrut[b]}}
\newcommand*\circled[1]{\tikz[baseline=(char.base)]{
		\node[shape=circle,draw,inner sep=2pt,fill=white] (char) {#1};}}
\newcommand*\circledText[1]{\tikz[baseline=(char.base)]{
		\node[shape=circle,draw,inner sep=2pt,fill=white] (char) {\hyperref[fig:WorkflowDiagram]{#1}};}}

\providecommand{\keywords}[1]{\textbf{\textit{Keywords ---}} #1}

\providecommand{\acknowledgements}[1]{\textbf{\textit{Acknowledgements ---}} #1}

\crefname{secname}{Section}{Section}
\Crefname{secname}{Sec.}{Sec.}
\crefname{section}{Section}{Sections}
\crefname{subsection}{Subsection}{Subsections}

\author{Johannes Rettberg\thanks{Institute of Engineering and Computational Mechanics, University of Stuttgart, Pfaffenwaldring 9, 70569 Stuttgart, Germany. (\url{johannes.rettberg,alexander.brauchler,pascal.ziegler,joerg.fehr@itm.uni-stuttgart.de})} \and Dominik Wittwar \thanks{Institute of Applied Analysis and Numerical Simulation, University of Stuttgart, Pfaffenwaldring 57, 70569 Stuttgart, Germany. (\url{dominik.wittwar,patrick.buchfink,haasdonk@mathematik.uni-stuttgart.de})} \and Patrick Buchfink\footnotemark[2] \and Alexander Brauchler\footnotemark[1] \and Pascal Ziegler\footnotemark[1] \and Jörg Fehr\footnotemark[1] \and Bernard Haasdonk\footnotemark[2]}
\title{Port-Hamiltonian Fluid-Structure Interaction Modeling and Structure-Preserving Model Order Reduction of a Classical Guitar}
\date{March 11, 2022}

\begin{document}
\maketitle
\begin{abstract}
\small
\textbf{Abstract} A fluid-structure interaction model in a port-Hamiltonian representation is derived for a classical guitar. We combine the laws of continuum mechanics for solids and fluids within a unified port-Hamiltonian (pH) modeling approach by adapting the discretized equations on second-order level in order to obtain a damped multi-physics model. The high-dimensionality of the resulting system is reduced by model order reduction. The article focuses on pH-systems in different state transformations, a variety of basis generation techniques as well as structure-preserving model order reduction approaches that are independent from the projection basis. As main contribution a thorough comparison of these method combinations is conducted. In contrast to typical frequency-based simulations in acoustics, transient time simulations of the system are presented. The approach is embedded into a straightforward workflow of sophisticated commercial software modeling and flexible in-house software for multi-physics coupling and model order reduction.
\end{abstract}


%
%
\section{Introduction}
\label{ch:Introduction}
The modeling and simulation of complex technical systems is indispensable for developing and operating new products. Different physical domains are suitably combined to create entirely new functionalities and obtain a more realistic model and a better understanding of the overall system. Considering multiple physical phenomena poses additional challenges in the system description as these systems need to be reasonably coupled. Their interaction must be taken into account for meaningful simulation results and a deeper understanding of the underlying system behavior. The port-Hamiltonian (pH) representation is based on the idea of interconnecting Hamiltonian systems and including dissipation in the formulation, for this reason the energy is used as the \textit{lingua franca} at the coupling ports \cite{SchaftJeltsema14}. This approach is attracting increasing attention due to its ability to describe complex multi-physics systems in the context of automated modeling and provides the basis for a unified description of modeling, analysis, control, optimization, and model order reduction. Furthermore, under mild assumptions pH-systems implicitly ensure important system properties
such as passivity, stability and energy conservation \cite{ChaturantabutEtAl16,RashadEtAl20}.

Fluid-structure interaction (FSI) characterizes the interaction of a moving or deformable mechanical structure with an internal or surrounding fluid \cite{BungartzSchaefer06}. The consideration of FSI plays a crucial role in different applications, e.g. the bending of aircraft wings, the blood flow in large veins or the lubricant flow in ball-bearings \cite{Richter17}. In the current study we want to analyze an acoustic system in which FSI is also of great importance for its proper technical function. 

\begin{figure}
	\centering
	\input{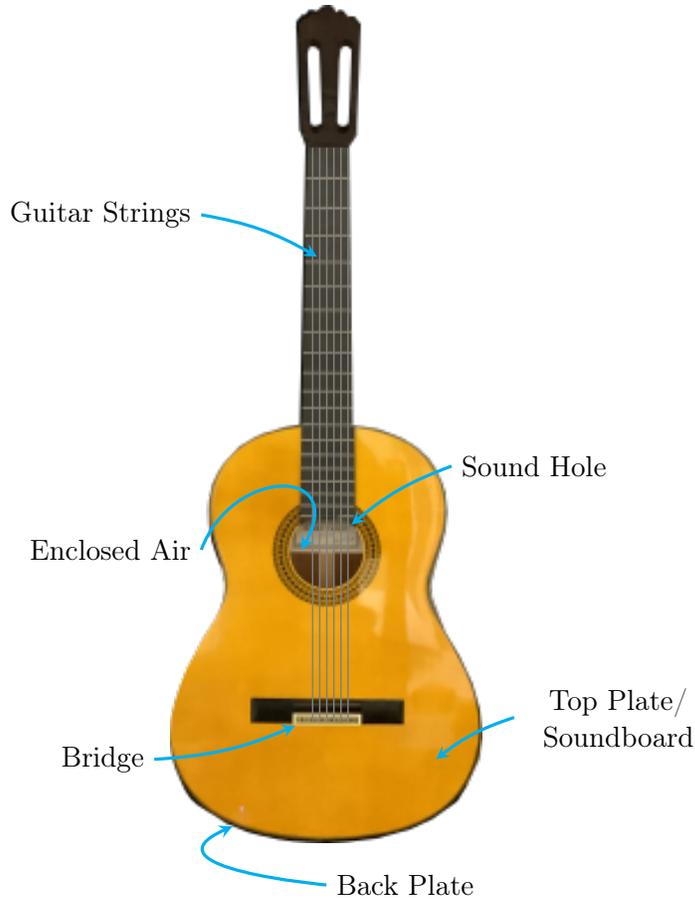}
	\caption{classical guitar}
	\label{fig:AcousticGuitar}
\end{figure}

A classical guitar is a popular instrument consisting of various components whose most important ones are briefly described in order to understand the physical context of the mathematical problem, see Fig.~\ref{fig:AcousticGuitar}. The musician plucks the guitar strings which thereupon start to vibrate. Those vibrations are transferred through the bridge into the soundboard which is the top plate of the guitar body. The guitar body coupled with the enclosed air acts as a resonator for the strings' vibration. The guitar's composition and material properties influence the volume and timbre of the air escaping through the sound hole.
The classical guitar represents an impressive example of a multi-physics system where the guitar body is strongly coupled with the enclosed air and the coupling effects cannot be neglected for a realistic representation of the behavior. The modeling of a guitar without strings can be found e.g. in \cite{BrauchlerEtAl21,Christensen82} and the modeling of strings  in \cite{BrauchlerEtAl20,DucceschiBilbao16}. Port-Hamiltonian representations of acoustic models are presented for instance in \cite{SilvaEtAl19} for the vocal apparatus and for a Rhodes piano in \cite{FalaizeHelie17}. The modeling procedure of this work follows closely the approach presented in \cite{BrauchlerEtAl21}.

In the current study this FSI problem is derived from the partial differential equation of linear elasticity for the mechanical structure \cite{TchonkovaSture01} and the wave equation for the fluid 
\cite{Blackstock00}. The system is spatially discretized with the finite element (FE) method in order to obtain a coupled second-order ordinary differential equation (ODE) system \cite{ZienkiewiczTaylorFox13}. Structural damping obtained from measurements is included into the model \cite{BrauchlerEtAl21}. From this basis, a port-Hamiltonian representation is attained by first transforming the system into a first-order system and then applying an adapted state transformation. The coupling effects of the multi-domain system are illustrated in a transient time simulation which gives further insight to classical frequency-based analyses in acoustics.

The spatial discretization via FE methods typically leads to high-dimensional systems that are ill-suited for multi-query simulations, e.g.\ for usage in optimization loops. Hence, model order reduction (MOR) methods are essential to approximate the system in a subspace of a much lower dimension for a faster and more efficient simulation \cite{Antoulas05}. One major novelty is the presentation of a comprehensive sensitivity analysis for a variety of combinations of projection methods and basis generation methods, i.e. the calculation of projection matrices, for reducing the FSI problem. The projection methods focus on basis-independent structure-preserving\footnote{Note that the term \textit{structure} has an ambiguity in this paper. On the one hand, it describes the mechanical structure of the wooden guitar plates as part of the fluid-structure interaction. On the other hand, it likewise describes the \textit{mathematical} structure of the pH-system, which should be preserved through the model reduction process. To avoid misunderstandings, the expression \textit{mechanical} structure is used in the context of the FSI.} reduction \cite{WolfEtAl10,EggerEtAl21} to preserve the important system properties in comparison to non-structure-preserving variants. The basis generation comprises approaches that are based on the system matrices, such as modal truncation \cite{Lehner07} or moment-matching via Krylov subspaces \cite{Freund03}, and data-based techniques based on time-response snapshots, e.g. Proper Orthogonal Decomposition (POD) \cite{Volkwein13} and Proper Symplectic Decomposition (PSD) \cite{PengMohseni16,BuchfinkEtAl19}. For more information on recent developments in structure-preserving MOR procedures for both Hamiltonian and port-Hamiltonian systems, we refer to \cite{BeattieEtAl21,HesthavenEtAl21,AltmannEtAl21,HesthavenPagliantini21} and the references therein.

\begin{figure}[h!]
	\centering
	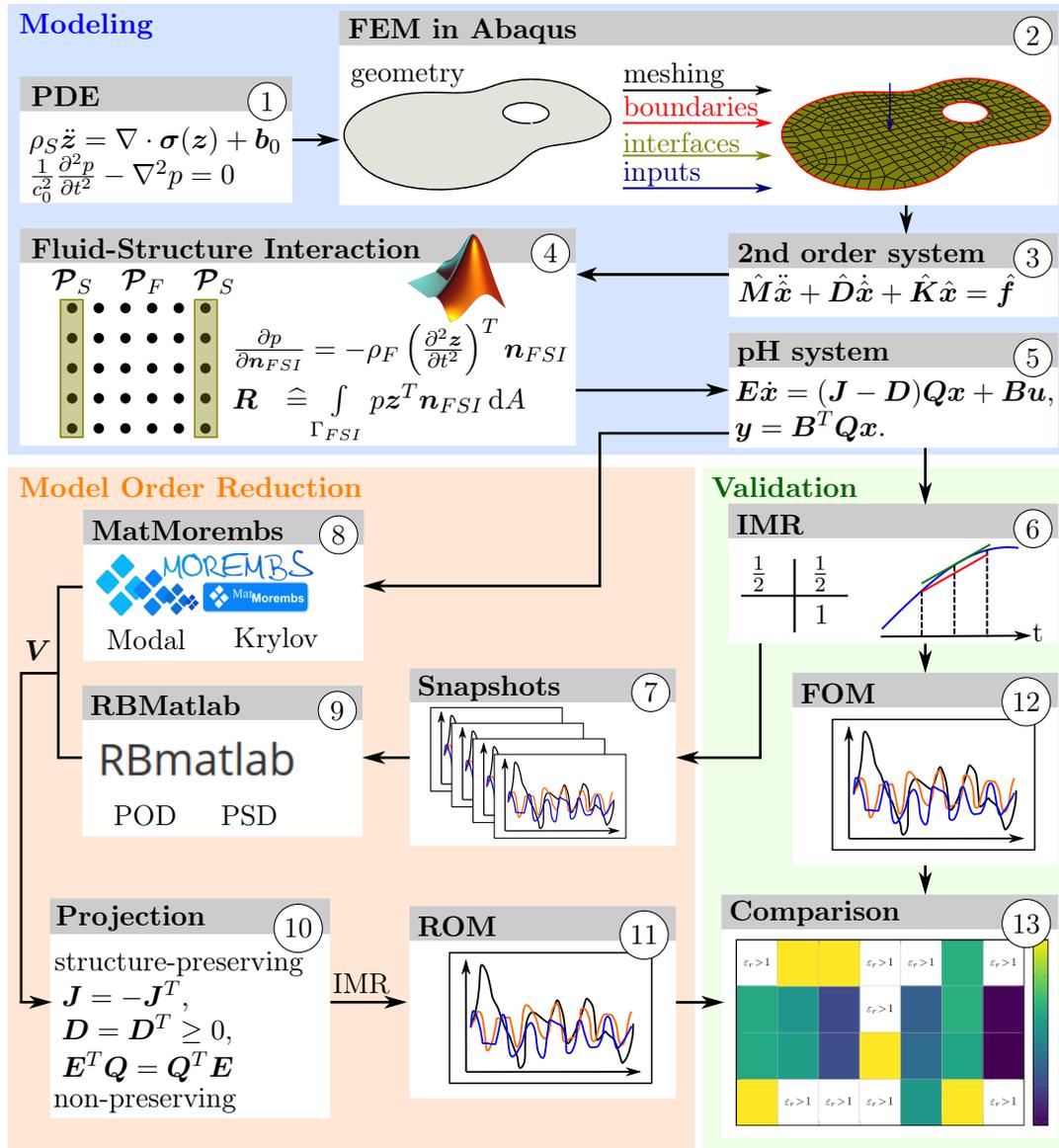
	\caption{workflow diagram}
	\label{fig:WorkflowDiagram}
\end{figure}

The paper is organized as follows. In Section~\ref{ch:ModelingGuitar}, the dynamic equations for a fluid-structure interaction problem are derived and adapted on second-order level. Additionally, practical aspects from the enhanced modeling process as an ensemble between the commercial software \textit{Abaqus} and in-house code in \textit{Matlab} are described. In Section~\ref{ch:portHamiltonianFormulation}, the necessary background on port-Hamiltonian systems is given and the FSI model of the guitar in a port-Hamiltonian input-output state representation is derived. Based on this, further pH formulations are deduced and the transient time simulation of the full system is performed. Various non-common projection methods, both structure-preserving and non-structure-preserving, are presented in Section~\ref{ch:MORProjectionMethod}. In Section~\ref{ch:MORBasisGeneration}, we recall classical and introduce modern basis generation methods with adapted snapshot matrices and non-orthogonal bases. The comprehensive sensitivity analysis is given in Section~\ref{ch:Results}. Finally, we conclude the paper. An overview of the entire workflow is given in Fig.~\ref{fig:WorkflowDiagram}. In the remainder this workflow diagram will be referenced by a circled number which highlights the workflow area the current section is referring to, e.g. \circledText{1}.

%
%
\section{Modeling of a classical guitar}
\label{ch:ModelingGuitar}
When producing a tone with a guitar, many different guitar components interact with each other, with the interaction of the guitar body and enclosed air being of particular importance. An enhanced three mass model focuses on this interaction and serves as a comprehensive system on which the investigated methodologies are presented. The present chapter first describes the mathematical background of this fluid-structure interaction and then also addresses the practical aspects of its implementation using \textit{Abaqus} and \textit{Matlab}.

%
%
\subsection{Derivation of the dynamic equations for the fluid-structure interaction}
\label{ch:DerivationFSI}
The vibration of the \textbf{mechanical structure} can be described as an elastic motion of the wooden plates of the guitar that can be mathematically explained by the methods of continuum mechanics, cf.~\circledText{1}. The motion of a body is described by the motion of all material points $\bm{\mathcal{P}}_S$ for this body, where the current position of each point can be classified by its position $\bm R_0$ in the initial configuration and the displacement $\bm z(\bm R_0,t)$. Assuming small strains leads to the elastodynamic problem for the time interval ${I_t=[0,t_{\text{end}}]}$
\begin{subequations}
	\begin{alignat}{5}
	\text{force equilibrium:}&\quad\rho_{S}\ddot{\bm z}&&=\nabla\cdot\bm \sigma(\bm z)+\bm b_0,\label{eq:forceEqu} \quad &&\text{in}~I_t\times\bm{\mathcal{P}}_S\\	
	\text{boundary condition:}&\qquad\bm z&&=\bm z_b,\quad &&\text{in}~I_t\times\Gamma_z \label{eq:dispBound}\\	
	\text{coupling term:}&\quad\bm \sigma \bm n&&= \bm t ,\quad &&\text{in}~I_t\times\Gamma_{FSI} \label{eq:stressBound}\\
	\text{strain-displacement relation:}&\qquad\bm \varepsilon_s&&=\frac{1}{2}\left(\nabla \bm z + (\nabla \bm z)^T\right),\quad &&\text{in}~\bm{\mathcal{P}}_S \label{eq:strainDispRel}\\
	\text{constitutive equation (Hooke's law):}&\qquad\bm \sigma &&=\bm C_{\text{EM}} \bm \varepsilon_s,\quad &&\text{in}~\bm{\mathcal{P}}_S \label{eq:constEqu} \\
	\text{initial condition:}&\qquad\bm z&&=\bm z_0,\quad &&\text{in}~\bm{\mathcal{P}}_S \label{eq:initCondDisp}\\
	&\qquad\dot{\bm z}&&=\dot{\bm z}_0,\quad &&\text{in}~\bm{\mathcal{P}}_S \label{eq:initCondVel}
	\end{alignat}
\end{subequations}
known as the strong formulation. The fundamental laws and material-dependence are considered in the force equilibrium \eqref{eq:forceEqu} and the constitutive equation \eqref{eq:constEqu}. The density of the material $\rho_S$, the second temporal derivative of the displacements $\ddot{\bm z}$, the divergence of the stress tensor $\bm \sigma$ and the volume forces $\bm b_0$ form the force equilibrium \eqref{eq:forceEqu}, where $\nabla$ describes the nabla operator. The constitutive equation \eqref{eq:constEqu} involves the engineering strain $\bm \varepsilon_s$ and elastic modulus $\bm C_{\text{EM}}$. We restrict ourselves to the case of orthotropic, linear elastic material behavior which can be parametrized with of nine parameters: three Young's moduli, three Poisson ratios and three shear moduli to consider the behavior in each spatial direction \cite{ZienkiewiczTaylorZhu05}. $\Gamma_z$ describe the parts of the body surface where Dirichlet boundary conditions for the displacement $\bm z_b$ apply. Likewise $\Gamma_{FSI}$ describes the contact surface of the mechanical structure and the encased air on which the stresses are given by the surface traction $\bm t = p \bm{n}_{FSI}$ whose magnitude scales proportionally with the pressure $p$ of the air at the contact surface. Here $\bm n_{FSI}$ denotes the normal vector of unit length pointing into the direction of the encased air. The initial conditions at time $t=0$ are given for the displacements $\bm z_0$ and velocities $\dot{\bm{z}}_0$ \cite{TchonkovaSture01,TkachukBischoff13,ZienkiewiczTaylorFox13}.

For a numerical simulation the continuum is spatially discretized with the finite element method (FEM)
resulting in the symmetric mass matrix $\bm M_S\in\mathbb{R}^{N_S \times N_S}$ and the symmetric stiffness matrix $\bm K_S\in\mathbb{R}^{N_S \times N_S}$.
To model dissipation effects, i.e.\ the damping parameters, we employ the so-called Rayleigh damping
\begin{equation}
\bm D_S = \alpha\bm M_S + \beta \bm K_S
\end{equation}
with the coefficients $\alpha$ and $\beta$. Measurements of the damping ratio were conducted on a real guitar at different frequencies $\omega_D$, see \cite{BrauchlerEtAl21}, and the overdetermined equation system is solved by using a least square fit for the coefficients that results in $\alpha = 15.958~\frac{1}{\text{s}}$ and $\beta = 2.821\cdot10^{-6}~\text{s}$.

This complements the dynamic equation of the mechanical structure resulting in the second-order space discretized system, c.f.~\circledText{3}
\begin{equation}
\label{eq:DynEqStructure}
\bm M_S\ddot{\bm z} + \bm D_S \dot{\bm{z}} + \bm K_S\bm z = \bm f_S + \bm f_p
\end{equation}
with the damping matrix $\bm D_S\in\mathbb{R}^{N_S \times N_S}$ and the excitation forces $\bm f_S, \bm f_p \in\mathbb{R}^{N_S}$ for the mechanical structure, whereas the latter is the specific force stemming from the coupling term \eqref{eq:stressBound}.

\begin{figure}[htb]
	\centering
	\input{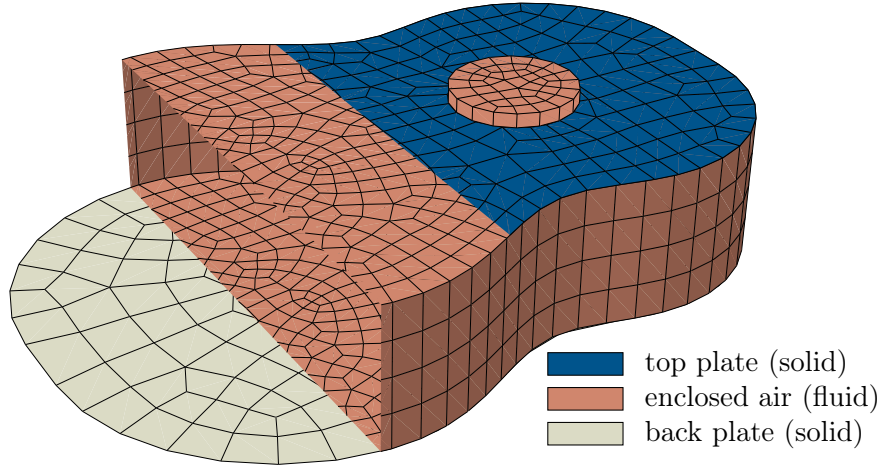}
	\caption{sectional view of the FEM multi-physics model of a classical guitar}
	\label{fig:GuitarCutView}
\end{figure}

The geometry and division of the guitar into the different domains can be seen in Fig.~\ref{fig:GuitarCutView}. The model consists of three different parts: the top plate or soundboard with a sound hole, the back plate and the air inside the guitar body in the shape and of the size of a classical guitar. Modeling the guitar behavior with these parts (in lumped parameter form) is known as a three mass model \cite{Christensen82}. The approach of \cite{BrauchlerEtAl21} is used as a reference to create an enhanced three mass model, which will be discussed further in Section~\ref{ch:ModelingAbaqusMatlab}.

The \textbf{fluid behavior}, c.f.~\circledText{1}, can be described with
\begin{subequations}
	\begin{alignat}{6}
	\text{acoustic wave equation:}&\quad \frac{1}{c_0^2}\frac{\partial^2p}{\partial t^2}-\nabla^2p &&=0,\quad &&\text{in}~I_t\times\bm{\mathcal{P}}_F 	\label{eq:AcousticWaveEq} \\
	\text{boundary conditions:}&\hspace{6.5em} p &&= p_b, \quad &&\text{in}~I_t \times \Gamma_b \\
	& \hspace{5.5em} \frac{\partial p}{\partial \bm{n}} && = 0, \qquad && \text{in}~ I_t \times \Gamma_n \\
	\text{coupling term:}&\hspace{4em} \frac{\partial p}{\partial \bm{n}_{FSI}} &&= -\rho_F \left(\frac{\partial^2 \bm{z}}{\partial t^2}\right)^T \bm{n}_{FSI}, \quad &&\text{in}~I_t \times \Gamma_{FSI} \label{eq:FluidCouplingTerm}
	\end{alignat}
\end{subequations}
where the acoustic wave equation is derived from the Euler equations for ideal, compressible fluids with Dirichlet boundary conditions $p_b~\text{on}~\Gamma_p$, zero Neumann boundary conditions, which model natural boundary conditions at the parts where the fluid is not encased by the structure, i.e. $\Gamma_n$ and non-zero Neumann boundary conditions at the interface, which describe the coupling relationship between the fluid and the mechanical structure. The force applied to the mechanical structure along the normal direction, i.e. $\frac{\partial p}{\partial \bm{n}_{FSI}}$ is equal in magnitude but of opposite direction compared to the force applied by the mechanical structure along the normal direction, i.e. $-\rho_F \left(\frac{\partial^2 \bm{z}}{\partial t^2}\right)^T \bm{n}_{FSI}$. The acoustic wave equation \eqref{eq:AcousticWaveEq} contains the speed of sound $c_0$ and the pressure $p$ on the domain $\bm{\mathcal{P}}_F$ \cite{Blackstock00,LerchEtAl09}. 

The wave equation is spatially discretized by using the Galerkin method with linear shape functions. This leads, c.f.~\circledText{3}, to the dynamic equations of the fluid in matrix form  \cite{HowardCazzolato17}
\begin{equation}
\label{eq:DynEqFluid}
\bm M_F\ddot{\bm p} + \bm K_F\bm p = \bm f_F + \bm f_{\bm{z}}
\end{equation}
with the symmetric mass matrix $\bm M_F\in\mathbb{R}^{N_F \times N_F}$, the symmetric stiffness matrix $\bm K_F\in\mathbb{R}^{N_F\times N_F}$, the excitation vector for the fluid $\bm f_F\in\mathbb{R}^{N_F}$ and the excitation vector $\bm f_{\bm{z}} \in \mathbb{R}^{N_F}$ stemming from the interaction with the mechanical structure, i.e.\ the coupling term \eqref{eq:FluidCouplingTerm}.

FSI is an area that describes the multi-physics coupling between the domains of fluid dynamics and structural mechanics. The coupling is considered a strong FSI if the motion of the solid influences the fluid and vice versa \cite{LerchEtAl09}. The classical guitar as the studied system allows the investigation of these strong couplings. In our case this coupling is expressed via the coupling terms \eqref{eq:stressBound} and \eqref{eq:FluidCouplingTerm}.

The effects of the mechanical structure on the fluid and vice versa, c.f.~\circledText{4}, result in the external forces $\bm f_{p}$ and $\bm f_{\bm{z}}$, respectively.  The former relates to the quantity
\begin{align}
\bm f_p \quad \widehat{=} \int\limits_{\Gamma_{FSI}} \bm {t}^T  \bm{z}_{\mathrm{test}} \, \mathrm{d} A = \int\limits_{\Gamma_{FSI}}  p  \bm{z}_{\mathrm{test}}^T\bm{n}_{FSI} \,\mathrm{d} A
\end{align}
whereas the latter corresponds to the quantity
\begin{align}
\bm f_{\bm z} \quad \widehat{=} \int\limits_{\Gamma_{FSI}} - \rho_F \left(\frac{\partial^2 \bm{z}}{\partial t^2}\right)^T \bm{n}_{FSI} p_{\mathrm{test}} \, \mathrm{d}A.
\end{align}
Hence, both $\bm{f}_{p}$ and $\bm{f}_{\bm{z}}$  can be represented as a matrix vector product
\begin{align}
\bm f_{p} = \bm R \bm{p} \qquad \text{ and } \qquad \bm f_{\bm{z}} = -\rho_F \bm R^T \ddot{\bm{z}},
\end{align}
where $\bm R$ describes the coupling matrix \cite{HowardCazzolato17} representing the term
\begin{align}
\bm R \quad \widehat{=} \int\limits_{\Gamma_{FSI}} p_{\mathrm{test}}  \bm{z}_{\mathrm{test}}^T\bm{n}_{FSI} \, \mathrm{d}A.
\end{align}
However, in the present case of non-matching meshes, the nodal values of the acoustic and structural domains, which are part of the interface, are related to each other with linear interpolation.

Consequently, the separate equations for the mechanical structure \eqref{eq:DynEqStructure} and the fluid \eqref{eq:DynEqFluid} can be combined into the following coupled system:

\begin{equation}
\label{eq:FSIMatrixEqUnsym}
\underbrace{\begin{bmatrix}
	\bm M_S&\bm 0\\\rho_F\bm R^T&\bm M_F
	\end{bmatrix}}_{\hat{\bm M}\in\mathbb{R}^{\hat{N}\times \hat{N}}}
\underbrace{\begin{bmatrix}
	\ddot{\bm z}\\\ddot{\bm p}
	\end{bmatrix}}_{\ddot{\hat{\bm x}}\in\mathbb{R}^{\hat{N}}}+
\underbrace{\begin{bmatrix}
	\bm D_S&\bm 0\\\bm 0&\bm 0
	\end{bmatrix}}_{\hat{\bm D}\in\mathbb{R}^{\hat{N}\times \hat{N}}}
\underbrace{\begin{bmatrix}
	\dot{\bm z}\\\dot{\bm p}\end{bmatrix}}_{\dot{\hat{\bm x}}\in\mathbb{R}^{\hat{N}}}
+
\underbrace{\begin{bmatrix}
	\bm K_S&-\bm R\\\bm 0&\bm K_F
	\end{bmatrix}}_{\hat{\bm K}\in\mathbb{R}^{\hat{N}\times \hat{N}}}
\underbrace{\begin{bmatrix}
	\bm z\\\bm p
	\end{bmatrix}}_{\hat{\bm x}\in\mathbb{R}^{\hat{N}}}
=\underbrace{\begin{bmatrix}
	\bm f_S(t)\\\bm f_F(t)
	\end{bmatrix}}_{\hat{\bm f}\in\mathbb{R}^{\hat{N}}}.
\end{equation}

Note that the mass and stiffness matrices, $\hat{\bm M}$ and $\hat{\bm K}$, are unsymmetric which leads to numerical issues when solving large sparse linear systems and to problems for the derivation of a port-Hamiltonian formulation of the system. Therefore, an adapted formulation of the system is derived in the following.

Instead of using the $\bm z$-$\bm p$-formulation \eqref{eq:FSIMatrixEqUnsym} a state transformation 
\begin{equation}
\label{eq:pDqTrans}
\bm p = \dot{\bm q}
\end{equation}
is applied, where instead of the pressure $\bm p$, the time derivative of the variable $\bm q$ is used which is in close relationship to the velocity potential. This approach is frequently used in the literature to symmetrize the system matrices \cite{WalleEtAl17,Everstine81}. 

A minor adaptation is used to make the pH matrices accomplish the mandatory skew-symmetry property after the transformation from a first-order to a second-order system: Inserting \eqref{eq:pDqTrans} into \eqref{eq:FSIMatrixEqUnsym}, taking the integral of the second line of \eqref{eq:FSIMatrixEqUnsym} and dividing by $ \rho_F$ (instead of $-\rho_F$ which is typically considered in literature) leads to the formulation 
\begin{equation}
\label{eq:FSIMatrixEqSym}
\begin{bmatrix}
\bm M_S&\bm 0\\\ \bm 0&\frac{1}{\rho_F}\bm M_F
\end{bmatrix}
\begin{bmatrix}
\ddot{\bm z}\\\ddot{\bm q}
\end{bmatrix}+
\begin{bmatrix}
\bm D_S&\bm R \\ -\bm R^T&\bm 0
\end{bmatrix}
\begin{bmatrix}
\dot{\bm z}\\\dot{\bm q}\end{bmatrix}
+
\begin{bmatrix}
\bm K_S&\bm 0\\\bm 0&\frac{1}{\rho_F}\bm K_F
\end{bmatrix}
\begin{bmatrix}
\bm z\\ \bm q
\end{bmatrix}
=\begin{bmatrix}
\bm f_S(t)\\\bm g_F(t)
\end{bmatrix}
\end{equation}
with 
\begin{equation}
\bm g_F(t) = \frac{1}{\rho_F}\int_{0}^{t}\bm f_F(\tau)\,d\tau.
\end{equation}

%
%
\subsection{Practical modeling aspects}
\label{ch:ModelingAbaqusMatlab}
The FE discretization is realized with the commercial software \textit{Abaqus}, see \cite{Abaqus16}, which allows us to design the complex geometry of the guitar, define material parameters and boundary conditions, and automatically mesh the geometry with sophisticated meshing tools, see~\circledText{2}.

The main goal of the current study is not a very accurate prediction of the guitar's behavior in accordance with experimental data but to rather investigate the performance of various formulations on an illustrative example for FSI. For this reason, the guitar is not modeled in full detail but still includes the most important aspects.

\begin{figure}[htb]
	\centering
\begingroup%
  \makeatletter%
  \providecommand\color[2][]{%
    \errmessage{(Inkscape) Color is used for the text in Inkscape, but the package 'color.sty' is not loaded}%
    \renewcommand\color[2][]{}%
  }%
  \providecommand\transparent[1]{%
    \errmessage{(Inkscape) Transparency is used (non-zero) for the text in Inkscape, but the package 'transparent.sty' is not loaded}%
    \renewcommand\transparent[1]{}%
  }%
  \providecommand\rotatebox[2]{#2}%
  \newcommand*\fsize{\dimexpr\f@size pt\relax}%
  \newcommand*\lineheight[1]{\fontsize{\fsize}{#1\fsize}\selectfont}%
  \ifx\svgwidth\undefined%
    \setlength{\unitlength}{311.81091927bp}%
    \ifx\svgscale\undefined%
      \relax%
    \else%
      \setlength{\unitlength}{\unitlength * \real{\svgscale}}%
    \fi%
  \else%
    \setlength{\unitlength}{\svgwidth}%
  \fi%
  \global\let\svgwidth\undefined%
  \global\let\svgscale\undefined%
  \makeatother%
  \begin{picture}(1,0.41538129)%
    \lineheight{1}%
    \setlength\tabcolsep{0pt}%
    \put(0,0){\includegraphics[width=\unitlength,page=1]{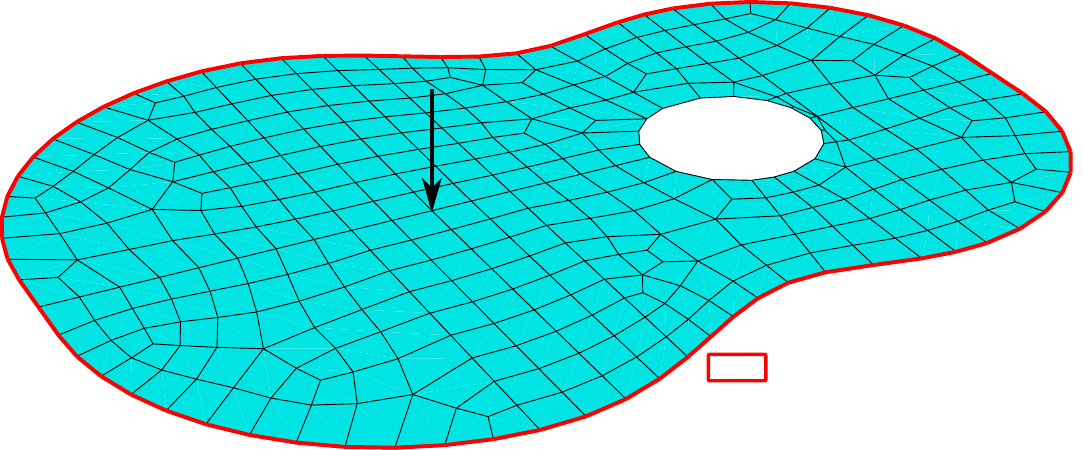}}%
    \put(0.71602337,0.06853746){\color[rgb]{0,0,0}\makebox(0,0)[lt]{\lineheight{1.25}\smash{\begin{tabular}[t]{l}boundary structure $\Gamma_z$\end{tabular}}}}%
    \put(0,0){\includegraphics[width=\unitlength,page=2]{SoundboardMeshEdit.pdf}}%
    \put(0.71697037,0.02394743){\color[rgb]{0,0,0}\makebox(0,0)[lt]{\lineheight{1.25}\smash{\begin{tabular}[t]{l}excitation force $\bm f_S$\end{tabular}}}}%
  \end{picture}%
\endgroup%

	\caption{meshed top plate}
	\label{fig:StructureMesh}
\end{figure}

\begin{figure}[htb]
	\centering
	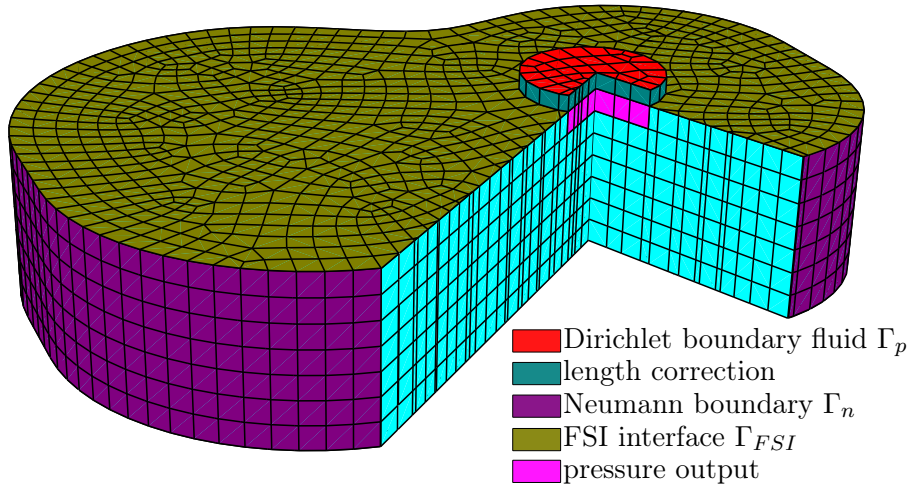
	\caption{meshed fluid}
	\label{fig:FluidMesh}
\end{figure}

The top and back plate are modeled with shell elements S4R5, see Fig.~\ref{fig:StructureMesh}. A summary of the finite element mesh properties used for the guitar model is presented in Tab.~\ref{tab:FEmeshProperties}. 
We use S4R5 element as it is accompanied by the advantage that the system size is reduced since there are less DOFs compared to the more common six DOFs S4 element. Furthermore, for the mechanical structure a six degree of freedom (DOF) formulation in \textit{Abaqus} would lead to a singular matrix $\bm M_S$ which would cause additional numerical issues and would lead to a differential-algebraic pH-system. 
The air model is enhanced by adding length correction to the fluid domain to account for effects of the surrounding air \cite{EzcurraEtAl05}, see Fig.~\ref{fig:FluidMesh}.

\begin{table}[htb]
	\centering
	\caption{finite element mesh properties}	
	{\begin{tabular}{lccc}
			\toprule
			&top plate&fluid&back plate \\ \midrule
			element type &shell S4R5&acoustic AC3D8& shell S4R5 \\
			nodes per element & 4 & 8 & 4 \\
			DOFs per node &5 &1 & 5 \\
			DOF types & 3 trans, 2 rot& 1 pressure& 3 trans, 2 rot \\ \midrule
			number of elements &343&2434&149 \\
			number of nodes &387&3226&171 \\
			number of DOFs & 1554 & 3170 & 900 \\ \bottomrule		
	\end{tabular}}
	\label{tab:FEmeshProperties}
\end{table}

As stated in the mathematical derivation of the mechanical structure and the fluid, some boundary conditions need to be considered that are in accordance with the physical behavior of the guitar. Therefore, the mechanical structure is constrained at the edges with homogeneous Dirichlet boundary conditions for the displacement variables. The rotational DOFs are not affected by this condition. This kind of boundary condition is considered as \textit{simply supported}. 
The fluid is mainly characterized by the coupling with the solid and the reflection from the sidewalls, due to homogeneous Neumann boundary conditions on $\Gamma_n$. However, the boundary condition at the top of the sound hole needs to be defined different since there is a free surface. This can be achieved by simply using the assumption of zero pressure \cite{ZienkiewiczTaylorZhu05}. The boundary conditions are marked in red and purple in Figs.~\ref{fig:StructureMesh} and \ref{fig:FluidMesh}.

The outcome of the modeling process in \textit{Abaqus} are the mass and stiffness matrices for the mechanical structure \eqref{eq:DynEqStructure} and the fluid \eqref{eq:DynEqFluid}. The further steps of modeling and simulation are performed in \textit{Matlab}. This change in the software tool is motivated by various advantages: (a) The modeling in \textit{Matlab} gives a better insight into the model especially due to the independent coupling procedure. In addition, (b) the obtained flexibility enables the system to be modified in the first place and thus to be reformulated into a pH-system. Furthermore, (c) this choice also gives us the freedom to decide on the integrator and the associated simulation properties.

In a real scenario, the guitarist plucks the string which introduces vibrations that are transferred through the bridge into the soundboard. Since the string and bridge are not modeled, the influence of these components is integrated as a force input at the position where the bridge would normally be located. A standard tuning of a guitar will create sounds in the range of approximately 82~Hz on the lowest pitch and 320~Hz at the highest pitch for the unfretted strings. The frequency range is chosen this way, since the objective is to demonstrate the methods on unfretted strings and the harmonics are neglected. For this reason, there will be a sine wave excitation force with
\begin{equation}
\label{eq:inputSineWave}
u = \hat{u}\sin(\omega t)
\end{equation} with $\omega=2\pi f$ in the frequency range $f=[82,320]~\text{Hz}$ and with an amplitude of $\hat{u}=1~\text{N}$ that acts on one node of the top plate.
As outputs, different important nodes of the mesh are considered. Those nodes include the aforementioned excitation node where the force acts on, compare Fig.~\ref{fig:StructureMesh}, as well as the volume-weighted integral of the pressure over the fluid nodes close to the sound hole of the guitar, compare Fig.~\ref{fig:FluidMesh}, and a central mechanical structure node of the back plate.

%
%
\section{Port-Hamiltonian formulation}
\label{ch:portHamiltonianFormulation}
Port-Hamiltonian (pH) systems were initially developed to describe a unified approach for systems from different physical domains that allow network modeling, see references in \cite{BreitenEtAl20,RashadEtAl20}. 
The pH formulation describes an energy-based network modeling approach with the energy as common quantity of various physical systems that interact through ports \cite{SchaftJeltsema14}. 

The pH structure can imply the underlying physical principles such as conservation laws \cite{BeattieEtAl18}. The input-output state form of a pH-system with neglected feed-through terms appears as 
\begin{equation}
\label{eq:pHSystemGeneral}
\begin{aligned}
\dot{\bm x}(t) &= (\bm J-\bm D)\nabla\mathcal{H}(\bm x(t)) + \bm B\bm u(t) \\
\bm y(t) &= \bm B^T\nabla\mathcal{H}(\bm x(t)) 
\end{aligned}
\end{equation}
where the function $\mathcal{H}(\bm x(t))$ describes the Hamiltonian as an energy function of the system. The Hamiltonian depends on the state vector $\bm x(t)\in\mathbb{R}^{N}$  with the system order $N = 2\hat{N}$. The matrix $\bm J\in\mathbb{R}^{N\times N}$ reflects the interconnection of the internal energy storage elements while the dissipation matrix $\bm D\in\mathbb{R}^{N\times N}$ describes the energy losses in the system. The matrix $\bm B\in\mathbb{R}^{N\times m}$ is the port or input matrix that specifies the manner in which energy enters or leaves the system. The Hamiltonian and the port matrix motivate the term pH-system. The pH-system matrices need to satisfy the requirement that $\bm J = -\bm J^T$ is skew-symmetric and $\bm D = \bm D^T\geq0$ is symmetric positive semidefinite. Furthermore, the variables $\bm u(t),\bm y(t)\in\mathbb{R}^{m}$ describe the time-dependent input and output vectors, respectively \cite{SchaftJeltsema14}.

Besides the modeling aspects, the pH formulations come with additional advantages as they implicitly incorporate important system properties. The pH-system \eqref{eq:pHSystemGeneral} is a generalization of a classical Hamiltonian system as the conservation of energy is generalized to a dissipation inequality
\begin{equation}
\label{eq:dissipationInequality}
\mathcal{H}(\bm x(t_1))-\mathcal{H}(\bm x(t_0))\leq\int_{t_0}^{t_1}\bm y(t)^T\bm u(t)\,dt\quad\text{with}\quad t_1>t_0
\end{equation}
which in combination with the assumption that the Hamiltonian is quadratic and strictly positive $\mathcal{H}(\bm x)>0$ for all $\bm x$ leads to the properties that the pH-system is both passive and stable \cite{ChaturantabutEtAl16}. Note, that the input and output entries of \eqref{eq:pHSystemGeneral} belong to the same variables defined by the port matrix $\bm B$. One can choose arbitrary quantities of interest (QoI) as outputs with appropriately chosen matrices but these output quantities may not satisfy the dissipation inequality \eqref{eq:dissipationInequality} and only serve the purpose of data observation.

Depending on the specific application and the coordinate system for which the problem is formulated, a generalized linear port-Hamiltonian system in the so-called descriptor form \cite{BeattieEtAl18},
\begin{equation}
\label{eq:pHDescriptorSystem}
\begin{aligned}
\bm E\dot{\bm x} &= (\bm J - \bm D) \bm Q \bm x+ \bm B \bm u, \\
\bm y &= \bm B^T\bm Q \bm x.
\end{aligned}
\end{equation}
with constant matrices $\bm E\in\mathbb{R}^{N\times N}$ and $\bm Q\in\mathbb{R}^{N \times N}$ that satisfy the condition $\bm E^T\bm Q = \bm Q^T\bm E$  in addition to the previously mentioned $\bm J$, $\bm D$ and $\bm B$ may also be considered.

In this case, the Hamiltonian of the system can be described as a quadratic energy function
\begin{equation}
\mathcal{H}(\bm x) = \frac{1}{2} \bm x^T\bm E^T\bm Q\bm x
\end{equation}

For our setting, such a descriptor system is obtained by transforming the second-order system \eqref{eq:FSIMatrixEqSym} into first-order a pH-system \eqref{eq:pHDescriptorSystem} via the new state vector $\bm x =\begin{bmatrix} \bm z^T & \bm q^T & \dot{\bm z}^T & \dot{\bm q}^T \end{bmatrix}^T$, which leads to the pH-system matrices
\begin{equation}
\label{eq:pHMatrices}
\begin{aligned}
\bm E &= \begin{bmatrix}
\bm I_{N_S} & \bm 0 & \bm 0 & \bm 0 \\
\bm 0 & \bm I_{N_F} & \bm 0 & \bm 0 \\
\bm 0 & \bm 0 & \bm M_S & \bm 0 \\
\bm 0 & \bm 0 & \bm 0 & \bm M_F/\rho_F
\end{bmatrix}\\
\bm J &= \begin{bmatrix}
\bm 0 & \bm 0 & \bm I_{N_S} & \bm 0 \\
\bm 0 & \bm 0 & \bm 0 & \bm I_{N_F} \\
-\bm I_{N_S} & \bm 0 & \bm 0 & \bm R \\
\bm 0 & -\bm I_{N_F} & -\bm R^T & \bm 0
\end{bmatrix}\\
\bm B\bm u&= \begin{bmatrix}
\bm 0 & \bm 0 & \bm f^T_S & \bm g^T_F
\end{bmatrix}^T
\end{aligned}
\quad
\begin{aligned}
\bm Q &= \begin{bmatrix}
\bm K_S & \bm 0 & \bm 0 & \bm 0 \\
\bm 0 & \bm K_F & \bm 0 & \bm 0 \\
\bm 0 & \bm 0 & \bm I_{N_S} & \bm 0 \\
\bm 0 & \bm 0 & \bm 0 & \bm I_{N_F}
\end{bmatrix}\\
\bm D &= \begin{bmatrix}
\bm 0 & \bm 0 & \bm 0 & \bm 0 \\
\bm 0 & \bm 0 & \bm 0 & \bm 0 \\
\bm 0 & \bm 0 & \bm D_S & \bm 0 \\
\bm 0  & \bm 0 & \bm 0 & \bm 0
\end{bmatrix}\\
\bm x &= \begin{bmatrix}
\bm z^T & \bm q^T & \dot{\bm z}^T & \dot{\bm q}^T
\end{bmatrix}^T
\end{aligned}
\end{equation}
with $\bm E$ describing the kinetic energy components and $\bm Q$ containing potential energy terms, c.f.~\circledText{5}. The matrices $\bm I_{N_{F/S}}$ represent the identity matrix of size $N_{F/S}$. The pH matrices satisfy the mandatory properties
\begin{subequations}
	\label{eq:pHProperties}
	\begin{align}
	\quad \bm J &= -\bm J^T, \label{eq:P1} \tag{pH1}\\
	\quad \bm D &= \bm D^T\geq 0, \label{eq:P2} \tag{pH2}\\
	\quad \bm E^T\bm Q &= \bm Q^T\bm E. \label{eq:P3} \tag{pH3}
	\end{align}
\end{subequations}
The pH-system \eqref{eq:pHMatrices} will further be considered as the velocity formulation.
\subsection{Numerical simulation of full-order results}
\label{ch:ResultsTimeIntegration}
We discretize the system with respect to time, using the implicit midpoint rule (IMR), cf.~\circledText{6}, which is the Gauss-Legendre collocation method with 1 stage and of order 2 for ODEs. The IMR is a symplectic integrator and hence, preserves the symplectic structure and the quadratic Hamiltonian through the time integration \cite{MehrmannMorandin19}.

We use a step size of $\Delta t = 10^{-4}~\text{s}$ and simulate for $T = 0.1~\text{s}$ to include several oscillation periods for the whole frequency range. All calculations were carried out with \textit{Matlab}. An advantage of a fixed time step is the better comparability of the trajectories. Since the numerical error of an integrator depends on the step size, we avoid mixing the integration and model reduction error.

The simulation of a full-order model (FOM) of the guitar is important for different reasons. On the one hand, these results will give us an insight about the approximate guitar behavior and the physical plausibility of the modeling. On the other hand, the FOM results form the basis for the calculation of the approximation error in the energy norm. Additionally, the FOM results are required for the snapshot generation, cf.~\circledText{7}, for the data-based basis generation techniques.

\begin{figure}[htp]
	\centering
	\input{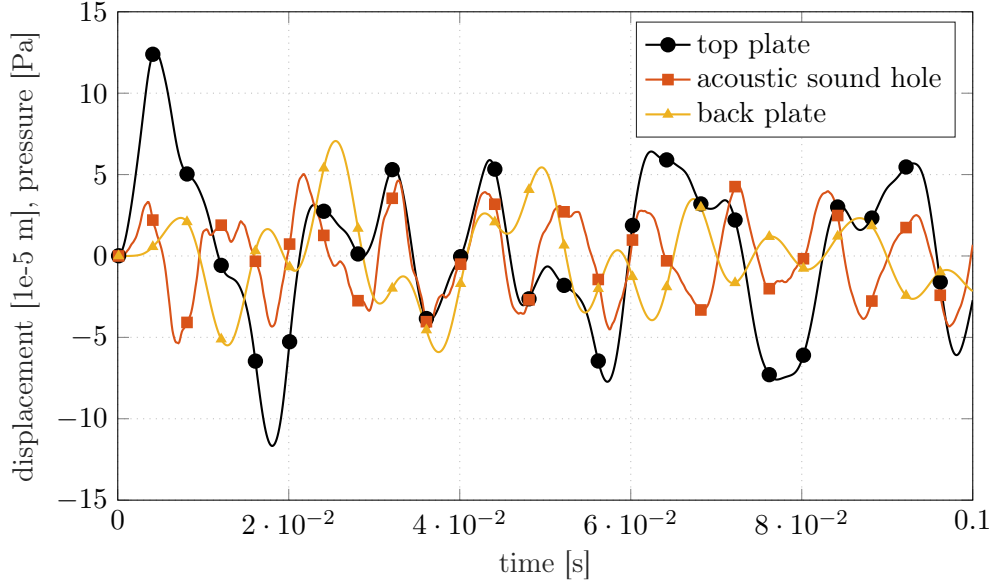}
	\caption{time simulation for one selected node per component at excitation frequency $f=100~\text{Hz}$}
	\label{fig:TimeSimulationM2T}
\end{figure}

We consider a parameter dependent excitation $\bm u = \bm u(\bm \mu)$ of the top plate that is modeled by a sinusoidal input where $\bm \mu \in [82,320]~\text{Hz}$ denotes the excitation frequency $f$. The time simulation of the FOM for an excitation frequency of $f=100~\text{Hz}$ is illustrated in Fig.~\ref{fig:TimeSimulationM2T}. It can be seen that as a result of the excitation, the top plate at the excitation node starts to vibrate at the same frequency since the guitar describes a linear time-invariant (LTI) system. The vibrations are transferred into the fluid whose pressure values start oscillating which corresponds to a sound of the guitar. The back plate  vibrates at a lower magnitude than the top plate which is justifiable by the transmission of energy through the fluid. It shall be emphasized that the back plate is only connected to the excited top plate via the fluid and hence, the strong coupling effects of the fluid can be demonstrated impressively.

%
%
\subsection{Reformulations of the pH-system}
\label{ch:SystemReformulations}
We derived our pH-FSI-system via a transformation of the second-order system and using the displacements, the integrated pressure and their derivatives as pH variables. In classical Hamiltonian systems the variables are defined in terms of the canonical position and momentum coordinates to account for additional symmetries \cite{ChaturantabutEtAl16}.  For this reason, a system in momentum formulation shall also be considered in the analysis. All of the following reformulations of the pH-system still satisfy the mandatory pH properties \eqref{eq:P1} to \eqref{eq:P3}. The coordinate transformation 
\begin{equation}
\bm x_{m}=\bm E\bm x
\end{equation}
yields such a pH-system in momentum formulation
\begin{equation}
\label{eq:pHDescriptorSystemMomentum}
\begin{aligned}
\dot{\bm x}_m &= (\bm J - \bm D) \bm Q \bm E^{-1}\bm x_m+ \bm B \bm u, \\
\bm y &= \bm B^T\bm Q \bm E^{-1} \bm x_m,
\end{aligned}
\end{equation}
where $\bm x_m$ contains the momentum instead of the velocities.

Since the PSD methods presented in Section~6 are based on a system formulation with canonical $\bm J$, i.e.
\begin{align*}
\bm{J} = \begin{bmatrix} \bm 0 & \bm{I}_{\hat{N}} \\ -\bm{I}_{\hat{N}} & \bm 0 \end{bmatrix},
\end{align*}
the canonical variants of the velocity formulation
\begin{equation}
\bm x_{c}=\bm P^{-1}\bm x
\end{equation}
and momentum formulation 
\begin{equation}
\bm x_{mc}=\bm P^{-1}\bm E\bm x
\end{equation}
will be considered in the sensitivity analysis, which yield the canonical system in velocity formulation
\begin{equation}
\label{eq:pHDescriptorSystemTrans}
\begin{aligned}
\bm P^T\bm E\bm P\dot{\bm x}_c &= \bm P^T(\bm J - \bm D)\bm P\bm P^{-1} \bm Q \bm P \bm x_c+ \bm P^T \bm B \bm u, \\
\bm y &= \bm B^T\bm Q \bm P \bm x_c
\end{aligned}
\end{equation}
and in momentum formulation
\begin{equation}
\label{eq:pHDescriptorSystemMomTrans}
\begin{aligned}
\bm P^T\bm P\dot{\bm x}_{mc} &= \bm P^T(\bm J - \bm D)\bm P\bm P^{-1} \bm Q \bm E^{-1} \bm P \bm x_{mc}+ \bm P^T \bm B \bm u, \\
\bm y &= \bm B^T\bm Q\bm E^{-1}\bm P \bm x_{mc}.
\end{aligned}
\end{equation}
The transformation matrices $\bm P,\bm P^{-1}$ depend on the coupling matrix $\bm R$ and are given by
\begin{equation}
\bm P = \begin{bmatrix}
\bm I_{N_S} & \bm 0 & \bm 0 & \bm R/2 \\
\bm 0 & \bm I_{N_F} & -\bm R^T/2 & \bm 0 \\
\bm 0 & \bm 0 & \bm I_{N_S} & \bm 0 \\
\bm 0 & \bm 0 & \bm 0 & \bm I_{N_F}
\end{bmatrix},
\qquad
\bm P^{-1} = \begin{bmatrix}
\bm I_{N_S} & \bm 0 & \bm 0 & -\bm R/2 \\
\bm 0 & \bm I_{N_F} & \bm R^T/2 & \bm 0 \\
\bm 0 & \bm 0 & \bm I_{N_S} & \bm 0 \\
\bm 0 & \bm 0 & \bm 0 & \bm I_{N_F}
\end{bmatrix}.
\end{equation}

The calculation of the FOM is computationally too expensive, especially if it comes to multi-query simulations. Therefore, MOR techniques are required and will be discussed in the following.

%
%
\section{Model order reduction - Projection methods}
\label{ch:MORProjectionMethod}
The size of a dynamical system model can be reduced by approximating its behavior in a subspace of a lower dimension. This model reduction comes at the cost of approximation errors. Projection-based reduction methods can be performed in various ways with different outcomes with respect to the preservation of different structural\footnote{Here, the term \textit{structure} refers to the \textit{mathematical} pH structure} pH properties, cf.~\circledText{10}.

All reduction methods used in our study are based on a projection-based reduction. The solution for $\bm x(t)$ is approximated in a subspace $\mathcal{V}$ of dimension $n\ll N$ which is described by a basis matrix $\bm V \in \mathbb{R}^{N\times n}$ with $\operatorname{colsp}(\bm V)=\mathcal{V}$ and leads to the approximation

\begin{equation}
\label{eq:ApproxVMatrix}
\bm x(t) \approx \bm V\bm x_r(t)
\end{equation}
with the reduced state vector $\bm x_r\in\mathbb{R}^n$.

Inserting \eqref{eq:ApproxVMatrix} into the pH-system \eqref{eq:pHDescriptorSystem} and using the Petrov-Galerkin condition for $\bm W\in\mathbb{R}^{N\times n}$ \cite{Antoulas05} yields the reduced system
\begin{equation}
\label{eq:pHReducedPetrovGal}
\begin{aligned}
\bm W^T\bm E\bm V\dot{\bm x_r}  &= \bm W^T(\bm J - \bm D) \bm Q \bm V \bm x_r + \bm W^T\bm B \bm u\\
\bm y &= \bm B^T\bm Q \bm V \bm x_r.
\end{aligned}
\end{equation}
of size $n \ll N$, where the matrix $\bm W$ determines the orthogonal projection direction.

A special case of the Petrov-Galerkin approach is the \textbf{Galerkin projection} with $\bm W = \bm V$ and takes the form
\begin{equation}
\label{eq:pHReducedGal}
\begin{aligned}
\bm V^T\bm E\bm V\dot{\bm x_r}  &= \bm V^T(\bm J - \bm D) \bm Q \bm V \bm x_r + \bm V^T\bm B \bm u\\
\bm y &= \bm B^T\bm Q \bm V \bm x_r.
\end{aligned}
\end{equation}
This projection does in general not preserve any of the underlying pH structure properties \eqref{eq:P1}-\eqref{eq:P3} of the system in our formulation.

The \textbf{quasi-Galerkin} projection is derived by inserting the term $\bm V\bm V^T$ under the assumption that $(\bm J-\bm D)\bm V\bm V^T\bm Q \approx (\bm J - \bm D)\bm Q$ to the standard Galerkin projection which yields
\begin{equation}
\label{eq:pHReducedQuasi}
\begin{aligned}
\underbrace{\bm V^T\bm{E}\bm V}_{\bm{E}_r}\dot{\bm x_r}  &= \bm V^T(\bm J - \bm D)\bm V\underbrace{\bm V^T \bm Q \bm V}_{\bm Q_r} \bm x_r + \bm V^T\bm B \bm u\\
\bm y &= \bm B^T\bm Q \bm V \bm x_r.
\end{aligned}
\end{equation}
This projection preserves \eqref{eq:P1} and \eqref{eq:P2} but it does not guarantee \eqref{eq:P3} and hence, the Hamiltonian of the reduced system $\bar{\mathcal{H}}(\bm x_r)\neq\mathcal{H}(\bm V\bm x_r)$ changes in the reduction.

Adapting the \textbf{pH structure-preserving} approach \cite{WolfEtAl10} with $\bm W = \bm Q\bm V$ to the case of a descriptor system leads to a reduced system
\begin{equation}
\label{eq:pHReducedpHpreserve}
\begin{aligned}
\bm V^T\bm Q^T\bm{E}\bm V\dot{\bm x_r}  &= \bm V^T\bm Q^T(\bm J - \bm D) \bm Q \bm V \bm x_r + \bm V^T\bm Q^T\bm B \bm u\\
\bm y &= \bm B^T\bm Q \bm V \bm x_r
\end{aligned}
\end{equation}
with the reduced matrices 
\begin{equation}
\begin{aligned}
\bm E_r &= \bm V^T\bm Q^T\bm E\bm V, & \bm J_r &= \bm V^T\bm Q^T\bm J\bm Q \bm V, \\
\bm D_r &= \bm V^T\bm Q^T\bm D\bm Q \bm V, & \bm Q_r &= \bm I_n.
\end{aligned}
\end{equation}
This preserves the pH properties \eqref{eq:P1}-\eqref{eq:P3} and does not change the Hamiltonian in the reduced system \cite{WolfEtAl10}.

The author in \cite{EggerEtAl21} presents a general framework for the numerical approximation of evolution problems which is also suitable for pH-systems and preserves the underlying Hamiltonian structure. We will call this approach \textbf{energy-stable}.
Applying this approach to a pH-system leads to a reduced system
\begin{equation}
\label{eq:pHReducedEgger}
\begin{aligned}
\bm V^T\bm{E}^T(\bm J-\bm D)^{-1}\bm{E}\bm V\dot{\bm x_r}  &= \bm V^T\bm{E}^T\bm Q \bm V \bm x_r + \bm V^T\bm{E}^T(\bm J-\bm D)^{-1}\bm B \bm u,\\
\bm y &= \bm B^T\bm Q \bm V \bm x_r.
\end{aligned}
\end{equation}
In our case the inverse term can be easily calculated with
\begin{equation}
(\bm J-\bm D)^{-1} = \bm T^T\bm P^T(\bm J-\bm D)\bm P\bm T\quad\text{with}\quad\bm T=\begin{bmatrix}
\bm 0 & \bm I_{\hat{N}}\\
\bm I_{\hat{N}} & \bm 0
\end{bmatrix}
\end{equation}
which allows for a computational efficient implementation.

In order to assess the quality of the different projection methods, a direct comparison with the \textbf{best approximation} in the considered subspace $\operatorname{colsp}(\bm V)$ should be examined. Here, we use the best approximation with respect to the energy norm $\lVert \cdot \rVert_{\mathcal{H}}^2$ induced by the energy inner product $\langle \cdot , \cdot \rangle_{\mathcal{H}}$
\begin{align}
\lVert\bm x\rVert_{\mathcal{H}}^2 := \bm x^T \bm H \bm x = \bm x^T\bm E^T\bm Q\bm x = \langle \bm x , \bm x \rangle_{\mathcal{H}},
\end{align}
where $\bm{H} = \bm{E}^T\bm{Q}$ denotes the so-called energy matrix which is positive definite, since both $\bm{E}^T$ and $\bm{Q}$ are positive definite and their product is symmetric due to the third pH-property \eqref{eq:P3}.
The best approximation of $\bm x$ in the subspace $\operatorname{colsp}(\bm V)$ is then given by the projection $\Pi_{\bm V,\mathcal{H}}$ of $\bm x$ onto $\operatorname{colsp}(\bm V)$ with respect to the energy inner product. In terms of the energy matrix this can be expressed via
\begin{align}
\Pi_{\bm V,\mathcal{H}} \bm x = \bm V \left( \bm V^T \bm H \bm V\right)^{-1} \bm V^T \bm H \bm x
\end{align}
and we have
\begin{align}
\label{eq:BestApproximation}
\lVert\bm x - \Pi_{\bm V,\mathcal{H}} \bm x \rVert_{\mathcal{H}} = \min\limits_{ \hat{\bm{x}} \in \operatorname{colsp}(\bm V)} \lVert\bm x - \hat{\bm{x}} \rVert_{\mathcal{H}}
\end{align}

%
%
\section{Model order reduction - Basis Generation}
\label{ch:MORBasisGeneration}
For all projection-based reduction methods presented in Section~\ref{ch:MORProjectionMethod}, a basis matrix $\bm V$ is required. In the following, different basis generation approaches are presented with the goal to keep the approximation error as small as possible. The bases are generated by using the extensive software packages MatMorembs\footnote{\href{https://www.itm.uni-stuttgart.de/software/morembs/software_morembs_matmorembs/}{https://www.itm.uni-stuttgart.de/software/morembs/software\textunderscore morembs\textunderscore matmorembs/}} \cite{FehrEtAl18b} and RBmatlab\footnote{\href{https://www.morepas.org/software/rbmatlab/}{https://www.morepas.org/software/rbmatlab/}}, cf.~\circledText{8} and \circledText{9}.

For the basis generation method based upon \textbf{modal truncation}, the homogeneous solution for \eqref{eq:pHDescriptorSystem} with $\bm u(t)=\bm 0$  is solved with the ansatz function
\begin{equation}
\bm x(t) = e^{\lambda_it}\bm \phi_i\quad\text{with}\quad\lambda_i\in\mathbb{C},\quad \bm \phi_i\in\mathbb{C}^{N}
\end{equation}
which leads to the generalized eigenvalue problem
\begin{equation}
\bm A\bm \phi_i = \lambda_i\bm{E}\bm \phi_i
\end{equation}
with $\bm A := (\bm J-\bm D)\bm Q$ for the pH-system. The eigenvectors $\bm \phi_i$ are also called coupled eigenmodes and describe the deformation of the mechanical structure and pressure distribution in the fluid at the dedicated eigenfrequencies. 
Modal truncation uses only the most important eigenvectors as the projection basis \cite{Antoulas05}. In the context of the guitar, the most important eigenvectors coincide with the eigenvectors that belong to the lowest eigenfrequencies since those are the crucial eigenfrequencies for the sound emission. Hence, the modal projection basis arises as
\begin{equation}
\bm V_{mod} = \begin{bmatrix}
\bm \phi_1 & \dots & \bm \phi_n
\end{bmatrix}
\end{equation}
with $n\ll N$. In the current study, the eigenmodes belonging to the lowest eigenfrequencies are coupled eigenmodes that allow for dynamics in the mechanical structure and fluid simultaneously. In general, $\bm V_{mod}$ is not real-valued. In this case, the matrix $\bm V_{mod}^T$ has to be replaced with $\bm V_{mod}^H$ for the projection methods outlined in Section~\ref{ch:MORProjectionMethod}, where the subscript $H$ represents the conjugate transpose of a complex-valued matrix. 

The goal of \textbf{Krylov-based reduction} is the approximation of the Laplace transformed transfer function of the pH-system 
\begin{equation}
\bm H(s) = \bm B^T\bm Q (s\bm E-(\bm J-\bm D)\bm Q)^{-1}\bm B
\end{equation}
with the complex variable $s\in\mathbb{C}$ which is typically evaluated on the imaginary axis $s=i\omega$ with the circular frequency $\omega = 2\pi f$ where $f$ denotes the excitation frequency. 

The transfer function around an expansion point $s_0$ can be described with a Taylor series whose coefficients are called moments in this context. The first $J_b$ moments of the full and reduced system can be matched by using the block Krylov subspace of a pH-system

\begin{equation}
\operatorname{colsp}(\bm V_{Kry}) = \mathcal{K}_{J_b} ((\bm A-s_0\bm E)^{-1}\bm E,(\bm A-s_0\bm E)^{-1}\bm B).
\end{equation}
The approach can be generalized to calculate for multiple expansion points \cite{Lehner07}. 

The direct calculation leads to numerical issues since additional vectors can become linearly dependent. That is why the block Arnoldi algorithm is used for the calculation of Krylov subspaces. This Arnoldi algorithm consists of a LU decomposition, Gram-Schmidt orthogonalization and the iterative calculation of Krylov directions \cite{Lehner07}. Sometimes the approach is also called tangential interpolation \cite{Castagnotto18}.

The \textbf{Proper Orthogonal Decomposition (POD)} approach uses a different idea for the projection basis generation which is based on snapshots instead of system matrices and can therefore even be used for nonlinear systems \cite{Volkwein13,Pinnau08}. 
The POD starts with a set of vectors $\hat{\bm x}_i\in\mathbb{R}^N$ with $i=1,\dots,m$ assembled column-wise into a matrix
\begin{equation}
\hat{\bm X} := \begin{bmatrix}
\hat{\bm x}_1 & \dots & \hat{\bm x}_m
\end{bmatrix}.
\end{equation}
The goal is to approximate the information contained in the snapshot matrix $\hat{\bm X}$ by a set of vectors $\hat{\bm u}_i\in\mathbb{R}^N$ with $i=1,\dots,n$, which can be expressed in terms of an optimization problem
\begin{equation}
\min\limits_{ \hat{\bm{u}}_1,\dots,\hat{\bm{u}}_n} \sum\limits_{i=1}^{k} \left\| \hat{\bm{x}}_i - \sum\limits_{j=1}^n (\hat{\bm{x}}_i^T\hat{\bm{u}}_j)\hat{\bm{u}}_j \right\|^2 \quad \text{subject to }\quad\hat{\bm u}_i^T\hat{\bm u}_j=\delta_{ij}
\end{equation}
where $\delta_{ij}$ describes the Kronecker delta and hence, the vectors $\hat{\bm u}_i$ form an orthonormal basis \cite{Volkwein13}.

The optimization problem can be solved by solving the eigenvalue problem
\begin{equation}
\label{eq:PODOptProbEVProbLeftEigenvector}
\hat{\bm X}\hat{\bm X}^T\hat{\bm u}_i=\lambda_i \hat{\bm u_i}
\end{equation}
with the scaled correlation matrix $\hat{\bm X}\hat{\bm X}^T\in\mathbb{R}^{N\times N}$ which is computationally inefficient to calculate if $m < N$ due to its size. Therefore, the method of snapshots solves the eigenvalue problem
\begin{equation}
\hat{\bm X}^T\hat{\bm X}\hat{\bm v}_i = \lambda_i\hat{\bm v}_i
\end{equation}
of size $m\times m$ with the eigenvectors $\hat{\bm v}_i$ which saves computational time for $m\ll N$ \cite{Volkwein13,Haasdonk17}. The problem \eqref{eq:PODOptProbEVProbLeftEigenvector} then is solved by
\begin{equation}
\hat{\bm u}_i = \frac{1}{\sqrt{\lambda_i}}\hat{\bm X}\hat{\bm v}_i.
\end{equation}

By only taking the first $n$ POD-modes, the basis matrix can be formed as
\begin{equation}
\bm V_{\text{POD}}(\hat{\bm X}):=\begin{bmatrix}
\hat{\bm u}_1 & \dots & \hat{\bm u}_n
\end{bmatrix}.
\end{equation}

In the present system, the snapshot matrix $\hat{\bm X}$ consists of state variables of the pH-system $\bm x(t,\bm \mu)$ with $\bm \mu \in [82,320]~\text{Hz}$. The discrete trajectories that form the snapshot matrix are calculated from the full system with time instances of $\Delta t = 10^{-4}~\text{s}$ for a time span of $T=0.1~\text{s}$. Using the \textbf{full state} vector $\bm x = \begin{bmatrix}
\bm z^T & \bm q^T & \dot{\bm z}^T & \dot{\bm q}^T
\end{bmatrix}^T$ for calculating, the POD basis will further be called as $\bm V_{\text{POD,State}}$. 

Another basis will be generated by only using the \textbf{displacements} in the snapshot matrix $\hat{\bm X}_{\text{Disp}}=\begin{bmatrix}
\hat{\bm x}_{\text{Disp},1}&\dots&\hat{\bm x}_{\text{Disp},m} \end{bmatrix}$ of the adapted state vector $\hat{\bm x}^T_{\text{Disp,i}} = \begin{bmatrix}
\bm z^T_i&\bm q^T_i
\end{bmatrix}$ and assembling the full projection matrix afterwards as
\begin{equation}
\bm V_{\text{POD,Disp}}(\hat{\bm X}_{\text{Disp}}) = \begin{bmatrix}
\bm V_{\text{POD}}(\hat{\bm X}_{\text{Disp}}) & \bm 0 \\
\bm 0 & \bm V_{\text{POD}}(\hat{\bm X}_{\text{Disp}})
\end{bmatrix}.
\end{equation}
Note that the above is a special case of the so-called tangent lift  as presented in \cite{PengMohseni16}, where only information in the state is considered. Hence, using the basis for the velocity component might be ill-suited.

A further approach will be the division of the state vector into the four individual parts and calculating the POD-basis for each \textbf{individual trajectory} and assembling the basis as the block-diagonal matrix
\begin{equation}
\bm V_{\text{POD,Individual}} = \begin{bmatrix}
\bm V_{\text{POD}}(\hat{\bm X}_{\bm z}) & \bm 0 & \bm 0 & \bm 0 \\
\bm 0 & \bm V_{\text{POD}}(\hat{\bm X}_{\bm q}) & \bm 0 & \bm 0 \\
\bm 0 & \bm 0 & \bm V_{\text{POD}}(\hat{\bm X}_{\dot{\bm z}}) & \bm 0 \\
\bm 0 & \bm 0 & \bm 0 & \bm V_{\text{POD}}(\hat{\bm X}_{\dot{\bm q}})
\end{bmatrix}.
\end{equation}
with the snapshot matrices $\hat{\bm X}_{\bm z}$, $\hat{\bm X}_{\bm q}$, $\hat{\bm X}_{\dot{\bm z}}$ and $\hat{\bm X}_{\dot{\bm q}}$ corresponding to the respective components of the state $\bm x$.

Symplectic MOR is structure-preserving MOR for Hamiltonian systems. One approach to derive a symplectic reduced-order basis is the data-driven \textbf{Proper Symplectic Decomposition (PSD)} which is closely related to the POD \cite{PengMohseni16}. 

Assuming a suitable basis, a symplectic form $\omega_{2\hat{n}}\colon\mathbb{R}^{2\hat{n}}\times \mathbb{R}^{2\hat{n}}\to\mathbb{R}$ takes the canonical form 
\begin{equation}
\omega_{2\hat{n}}(\bm v_1,\bm v_2) = \bm v^T_1\mathbb{J}_{2\hat{n}}\bm v_2\quad\forall\bm v_1,\bm v_2\in\mathbb{V} \quad \text{with}\quad \mathbb{J}_{2\hat{n}}:=\begin{bmatrix}
\bm 0_{\hat{n}} & \bm I_{\hat{n}} \\ -\bm I_{\hat{n}} & \bm 0_{\hat{n}}
\end{bmatrix}
\end{equation}
where $\mathbb{J}_{2\hat{n}}$ is the Poisson matrix with the identity matrix $\bm I_{\hat{n}} \in \mathbb{R}^{\hat{n} \times \hat{n}}$ and the matrix of all zeros $\bm 0_{\hat{n}}\in\mathbb{R}^{\hat{n} \times \hat{n}}$ \cite{Silva08}. We call $\bm V\in\mathbb{R}^{2\hat{N}\times2\hat{n}}$ a symplectic matrix if
\begin{equation}
\bm V^T\mathbb{J}_{2\hat{N}}\bm V=\mathbb{J}_{2\hat{n}}.
\end{equation}

If a Petrov-Galerkin projection is used with $\bm W^T =\bm V^+$ where 
\begin{equation}
\bm V^+ = \mathbb{J}_{2\hat{n}}^T\bm V^T\mathbb{J}_{2\hat{N}} \in \mathbb{R}^{2\hat{n}\times 2\hat{N}}
\end{equation}
denotes the so-called symplectic inverse\footnote{Please note, that this symplectic inverse should not be confused with the Moore-Penrose inverse, which sometimes is denoted with the same symbol.}, then the Hamiltonian structure and therefore its Hamiltonian is preserved for purely Hamiltonian systems, i.e. in the case of $\bm{D} = 0$ in equation \eqref{eq:pHSystemGeneral}\cite{PengMohseni16,BuchfinkEtAl19}.

In analogy to the POD, the minimization problem appears as
\begin{equation}
\label{eq:PSDOptimization}
\begin{aligned}
&\underset{\bm V\in\mathbb{R}^{2\hat{N}\times 2\hat{n}}}{\text{minimize}} & & \lVert (\bm I_{2\hat{N}}-\bm V\bm V^+)\hat{\bm X}\rVert_F^2\\
& \text{subject to} & & \bm V^T\mathbb{J}_{2\hat{N}}\bm V=\mathbb{J}_{2\hat{n}},
\end{aligned}
\end{equation}
where the constraint guarantees that the reduced-order basis (ROB) is symplectic \cite{PengMohseni16}. In contrast to POD, there is no explicit solution procedure for the PSD optimization problem \eqref{eq:PSDOptimization}. The \textbf{PSD Complex SVD} approach is the solution of the PSD in the subset of symplectic, orthonormal ROBs \cite{BuchfinkEtAl19}. It uses an adapted complex snapshot matrix 
\begin{equation}
\bm C_s = \begin{bmatrix}
\tilde{\bm x}_1^s +i\dot{\tilde{\bm x}}_1^s & \dots & \tilde{\bm x}_m^s +i\dot{\tilde{\bm x}}_m^s
\end{bmatrix}\in \mathbb{C}^{\hat{N}\times m}\quad
\tilde{\bm x}_j^s = \begin{bmatrix}
\bm z_j \\ \bm q_j
\end{bmatrix}\quad 1\leq j\leq m
\end{equation} 
with the imaginary unit $i$. The minimization problem requires an auxiliary complex matrix $\bm U_{\bm C_s}\in\mathbb{C}^{\hat{N}\times\hat{n}}$ and takes the form 
\begin{equation}
\label{eq:MinProblemPSDComplexSVD}
\begin{aligned}
&\underset{\bm U_{\bm C_s}\in\mathbb{C}^{2\hat{N}\times 2\hat{n}}}{\text{minimize}} & & \lVert \bm C_s -\bm U_{\bm C_s}(\bm U_{\bm C_s})^H\bm C_s\rVert_F^2\\
& \text{subject to} & & (\bm U_{\bm C_s})^H\bm U_{\bm C_s}=\bm I_{\hat{n}}
\end{aligned}
\end{equation}
and generates the real basis matrix
\begin{equation}
\bm V_{\text{CSVD}}:=\begin{bmatrix}
\widetilde{\bm E} & \mathbb{J}_{2\hat{N}}\widetilde{\bm E}
\end{bmatrix},\quad \widetilde{\bm E} := \begin{bmatrix}
\operatorname{Re}(\bm U_{\bm C_s})\\
\operatorname{Im}(\bm U_{\bm C_s})
\end{bmatrix}.
\end{equation}
The specialty of the PSD Complex SVD is the choice of the auxiliary complex matrix $\bm C_s$ and the computation of $\widetilde{\bm E}$ from $\bm U_{\bm C_s}$. The solution of \eqref{eq:MinProblemPSDComplexSVD} is the POD for complex matrices based on the left-singular vectors of $\bm C_s$ that can be computed with a complex version of the SVD \cite{PengMohseni16}.

Most existing basis generation techniques, e.g. Complex SVD, generate a symplectic, orthonormal ROB. In \cite{BuchfinkEtAl19} a new symplectic, non-orthogonal basis generation technique based on the so-called \textbf{SVD-like decomposition} is introduced. 

There exists a SVD-like decomposition of the snapshot matrix $\hat{\bm X}\in\mathbb{R}^{2\hat{N}\times m}$ 
\begin{equation}
\hat{\bm X} = \bm S_s\bm D_s\bm Q_s,\quad
\bm D_s=\begin{blockarray}{ccccc}
\mLabel{p}&\mLabel{q}&\mLabel{p}&\mLabel{m-2p-q} \\
\begin{block}{[cccc]c}
\bm \Sigma_s & \bm 0 & \bm 0 & \bm 0 & \mLabel{p} \topstrut\\
\bm 0 & \bm I & \bm 0 & \bm 0 & \mLabel{q} \\
\bm 0 & \bm 0 & \bm 0 & \bm 0 & \mLabel{\hat{N}-p-q} \\
\bm 0 & \bm 0 & \bm \Sigma_s & \bm 0 & \mLabel{p} \\
\bm 0 & \bm 0 & \bm 0 & \bm 0 & \mLabel{q} \\
\bm 0 & \bm 0 & \bm 0 & \bm 0 & \mLabel{\hat{N}-p-q} \botstrut\\
\end{block}
\end{blockarray}, \quad \bm \Sigma_s = \operatorname{diag}(\sigma_1,\dots,\sigma_p)\in\mathbb{R}^{p\times p}
\end{equation}
with a symplectic matrix $\bm S_s\in\mathbb{R}^{2\hat{N}\times 2\hat{N}}$, a sparse and potentially non-diagonal matrix $\bm D_s\in\mathbb{R}^{2\hat{N}\times m}$, an orthogonal matrix $\bm Q_s\in\mathbb{R}^{m\times m}$ and the symplectic singular values $\sigma_i$. The rank of the snapshot matrix is $\operatorname{rank}(\hat{\bm X})=2p+q$.

The Frobenius norm of the snapshot matrix can be rewritten as 
\begin{equation}
\lVert \hat{\bm X}\rVert^2_F =\operatorname{Tr}(\hat{\bm X}\hat{\bm X}^T) = \sum_{i=1}^{p+q}(w_i)^2,\quad w_i=\begin{cases}
\sigma_i\sqrt{\lVert \bm s_i\rVert_2^2+\lVert \bm s_{\hat{N}+i}\rVert^2_2}&1\leq i\leq p \\
\lVert\bm s_i\rVert_2, & p+1 \leq i \leq p+q
\end{cases}
\end{equation}
where $\bm s_i$ is the $i$-th column of $\bm S_s$ and $w_i$ is called the weighted symplectic singular value \cite{BuchfinkEtAl19}. 

The goal is to choose the $k$ indices $i\in\mathcal{I}_{\text{SVD}}=\{i_1,\dots,i_k\}\subset\{1,\dots,p+q\}$ which have large contributions $w_i$ to the Frobenius norm with
\begin{equation}
\mathcal{I}_{\text{SVD}}=\underset{\begin{subarray}{c}
	\mathcal{I}\subset\{1,\dots,p+q\} \\
	\lvert\mathcal{I}\rvert=k
	\end{subarray}}{\operatorname{argmax}}\left(\sum_{i\in\mathcal{I}} (w_i)^2\right)
\end{equation}
which yields a ROB $\bm V_{\text{SVD}}\in\mathbb{R}^{2\hat{N}\times 2k}$ with $k$
pairs of columns $\bm s_i\in\mathbb{R}^{2\hat{N}}$ from $\bm S_s$ leading to reduced model order of $2k=n$ and the basis matrix
\begin{equation}
\bm V_{\text{SVD}}=\begin{bmatrix}
\bm s_{i_1}&\dots&\bm s_{i_k}&\bm s_{\hat{N}+i_1}&\dots&\bm s_{\hat{N}+i_k}
\end{bmatrix}.
\end{equation}
The SVD-like decomposition is constructed by computing an eigendecomposition of $\hat{\bm X}^T\bm{J}\hat{\bm X}$, for which the imaginary and real part of eigenvectors corresponding to complex conjugate eigenvalues result in a pair of
symplectic vectors. Alternatively, the approach introduced in \cite{Xu05} can be used. However, we chose the former method, as the later one can be quite computationally demanding.

%
%
\section{Results}
\label{ch:Results}
The previous sections featured several alternatives for model order reduction of a classical guitar FSI problem in a pH formulation. The different components that were considered are the 4 system formulations in Sec.~\ref{ch:SystemReformulations}, the 4 projection methods from Sec.~\ref{ch:MORProjectionMethod} and the 7 basis generation methods in Sec.~\ref{ch:MORBasisGeneration}, cf. Fig.~\ref{fig:WorkflowDiagram}. These categories are extended by choosing the reduced system size $n$, studying 11 sizes between $n=12$ and $n=400$ and the number of trajectories in the snapshot matrix for the data-based methods, where 6 alternatives were investigated, see Tab.~\ref{tab:ErrorOverSnapshots}. All representatives of each category can be combined with all others, resulting in a number of about $4\cdot4\cdot7\cdot11\cdot6\approx7000$ possible combinations. All of these combinations have been computed. We will give a focused view onto the most important outcomes in the following sensitivity analysis.

\subsection{Sensitivity analysis of different reduced-order models}
\label{ch:SensitivityROM}
In order to compare the different methods for reducing the FOM to a reduced-order model (ROM), cf.~\circledText{11}, an error measure is needed. The energy norm of the Hamiltonian
\begin{equation}
\lVert\bm x(t)\rVert_{\mathcal{H}}^2 = \bm x(t)^T\bm E^T\bm Q\bm x(t)
\end{equation}
is used to describe the energy value of a trajectory at time step $t$. The error $\varepsilon(t)=\lVert\bm x(t)-\bm V\bm x_r(t)\rVert_{\mathcal{H}}$ between the full and reduced model is measured in the energy norm. In order to obtain a scalar and relative measure, the norm is integrated over the simulation time interval $I_t=[0,t_{\text{end}}]$ and related to the full model which yields
\begin{equation}
\varepsilon_r = \frac{\int_{0}^{t_{\text{end}}}\lVert\bm x(t)-\bm V\bm x_r(t)\rVert_{\mathcal{H}}\,dt}{\int_{0}^{t_{\text{end}}}\lVert\bm x(t)\rVert_{\mathcal{H}}\,dt}
\end{equation}
for the relative error measure. It is important to use the relative error which allows for a comparison between different system inputs. Five randomly chosen frequencies in the standard tuning range \eqref{eq:inputSineWave} and zero initial conditions were used as system parameters for all experiments that were carried out. The mean value of their relative error calculations is displayed in the subsequent results.

\begin{figure}[htp]
	\centering
%
%

\pgfplotsset{compat=newest}
\usetikzlibrary{plotmarks}
\usetikzlibrary{arrows.meta}
\usepgfplotslibrary{patchplots}
\usepgfplotslibrary{groupplots}
\pgfplotsset{compat=newest}
\usepgfplotslibrary{groupplots}

\definecolor{mycolor1}{rgb}{0.85000,0.32500,0.09800}%
\definecolor{mycolor2}{rgb}{0.92900,0.69400,0.12500}%
\definecolor{mycolor3}{rgb}{0.49400,0.18400,0.55600}%
\begin{tikzpicture}
\begin{groupplot}[group style={group size= 2 by 1,
group name =PHEnSt,
horizontal sep = 1.4 cm,
},
legend columns = 4,
]
\nextgroupplot[%
title = pH,
width=1.779in,
height=2.567in,
scale only axis,
xmin=1,
xmax=7,
xtick={1,2,3,4,5,6,7},
xticklabels={{Modal},{Krylov},{POD-State},{POD-Disp},{POD-Indiv},{C-SVD},{SVD-like}},
xticklabel style={rotate=90},
ymode=log,
ymin=0.001,
ymax=1.053237,
yminorticks=true,
ylabel style={font=\color{white!15!black}},
ylabel={Relative Error},
axis background/.style={fill=white},
xmajorgrids,
ymajorgrids,
yminorgrids,
grid style={dotted},
every axis plot/.append style={thick},
]

\addplot [color=black, mark size=3pt, mark=*, mark options={solid, fill=black, black}, forget plot]
table[row sep=crcr]{%
	1	0.0976531\\
	2	0.02693616\\
	3	0.05965605\\
	4	0.02219798\\
	5	0.02462164\\
	6	0.1267772\\
	7	0.01198542\\
};
\addplot [color=mycolor1, mark size=2.5pt, mark=square*, mark options={solid, fill=mycolor1, mycolor1}, forget plot]
table[row sep=crcr]{%
	1	0.0976531\\
	2	0.07517606\\
	3	0.01140042\\
	4	1\\
	5	0.02226144\\
	6	0.0905212899999999\\
	7	0.002821226\\
};
\addplot [color=mycolor2, mark size=2.2pt, mark=triangle*, mark options={solid, fill=mycolor2, mycolor2}, forget plot]
table[row sep=crcr]{%
	1	0.0976531\\
	2	0.07303923\\
	3	0.06604698\\
	4	1.053237\\
	5	0.3930059\\
	6	0.1257139\\
	7	0.02763844\\
};
\addplot [color=mycolor3, mark size=1.8pt, mark=diamond*, mark options={solid, fill=mycolor3, mycolor3}, forget plot]
table[row sep=crcr]{%
	1	0.0976531\\
	2	0.06005757\\
	3	0.01114028\\
	4	1\\
	5	0.02226161\\
	6	0.0905209300000001\\
	7	0.002821227\\
};

\nextgroupplot[%
legend to name = named,
title = Energy--Stable,
width=1.779in,
height=2.567in,
scale only axis,
xmin=1,
xmax=7,
xtick={1,2,3,4,5,6,7},
xticklabels={{Modal},{Krylov},{POD-State},{POD-Disp},{POD-Indiv},{C-SVD},{SVD-like}},
xticklabel style={rotate=90},
ymode=log,
ymin=0.001,
ymax=10,
yminorticks=true,
axis background/.style={fill=white},
xmajorgrids,
ymajorgrids,
yminorgrids,
grid style={dotted},
every axis plot/.append style={thick},
]

\addplot [color=black, mark size=3pt, mark=*, mark options={solid, fill=black, black}]
table[row sep=crcr]{%
	1	0.0976531\\
	2	0.02098306\\
	3	0.02628112\\
	4	0.02219798\\
	5	0.01466514\\
	6	0.0421542199999999\\
	7	0.00940123299999999\\
};
\addlegendentry{Vel}

\addplot [color=mycolor1, mark size=2.5pt, mark=square*, mark options={solid, fill=mycolor1, mycolor1}]
table[row sep=crcr]{%
	1	0.0976531\\
	2	0.0593407099999999\\
	3	0.01039362\\
	4	2.839854\\
	5	0.01638676\\
	6	0.1141867\\
	7	0.002967628\\
};
\addlegendentry{Mom}

\addplot [color=mycolor2, mark size=2.2pt, mark=triangle*, mark options={solid, fill=mycolor2, mycolor2}]
table[row sep=crcr]{%
	1	0.0976531\\
	2	0.07459076\\
	3	0.04321814\\
	4	1.371127\\
	5	0.3385551\\
	6	0.05411288\\
	7	0.02358035\\
};
\addlegendentry{Vel\_Canon}

\addplot [color=mycolor3, mark size=1.8pt, mark=diamond*, mark options={solid, fill=mycolor3, mycolor3}]
table[row sep=crcr]{%
	1	0.0976531\\
	2	0.0580762199999999\\
	3	0.01056518\\
	4	2.839854\\
	5	0.01638674\\
	6	0.1141868\\
	7	0.002967627\\
};
\addlegendentry{Mom\_Canon}

\end{groupplot}
\path (PHEnSt c1r1.south) -- node[shift={(0,-2.43)}]{\ref*{named}} (PHEnSt c2r1.south);
\end{tikzpicture}%
	\caption{relative error of different basis generation methods and system formulations for $n=120$}
	\label{fig:SysForm_PH_ES}
\end{figure}
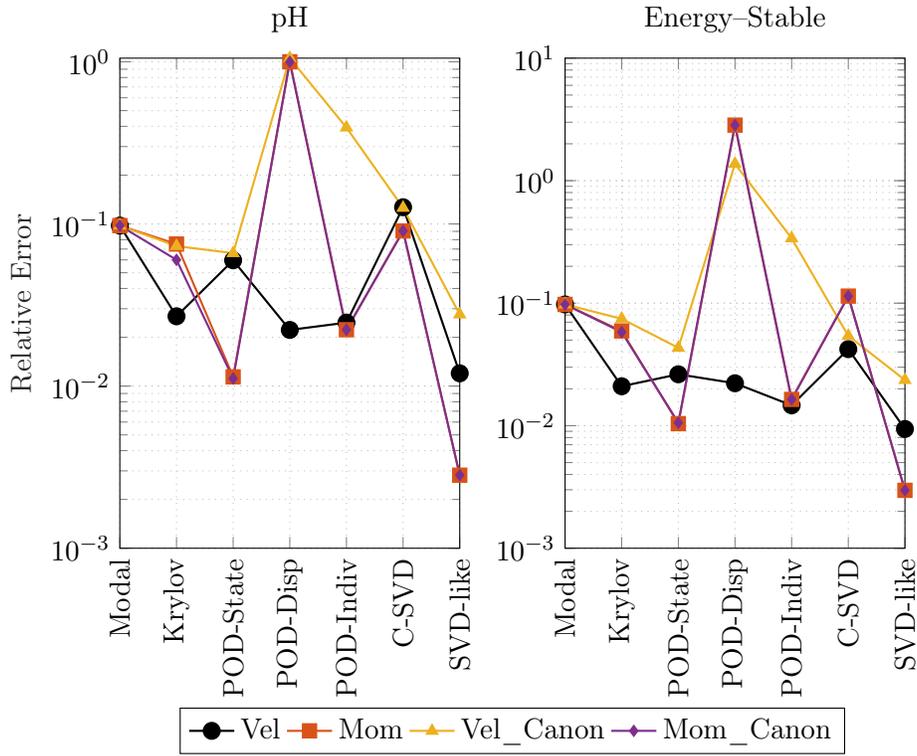

The relative error values for the four system formulations over all basis generation techniques and for the pH and energy-stable projections are illustrated in Fig.~\ref{fig:SysForm_PH_ES}. It can be seen that the displacement variant of POD is not a suitable option in the momentum formulation. This is due to the fact that the state derivatives are not considered in the snapshot matrix and hence, the inertia information gets lost. 
The canonical alternative for the velocity formulation \eqref{eq:pHDescriptorSystemTrans} generally leads to much higher errors and also in the momentum formulation the transformation \eqref{eq:pHDescriptorSystemMomTrans} is not increasing the performance. Furthermore, in the transformation process, the block structure of the pH matrices $\bm E$ and $\bm Q$ \eqref{eq:pHMatrices} gets lost. The momentum formulation for both, canonical and non-canonical, leads to the best overall results especially for POD-State and SVD-like. For these reasons, we will focus on the non-canonical momentum formulation in the remainder.

\begin{figure}[htp]
	\centering
%
%
\pgfplotsset{compat=newest}
\usetikzlibrary{plotmarks}
\usetikzlibrary{arrows.meta}
\usepgfplotslibrary{patchplots}

\definecolor{mycolor1}{rgb}{0.92900,0.69400,0.12500}%
\definecolor{mycolor2}{rgb}{0.46600,0.67400,0.18800}%
\definecolor{mycolor3}{rgb}{0.63500,0.07800,0.18400}%
%
\def\mystrut{\vphantom{hg}}
\pgfplotsset{
	legend image with text/.style={
		legend image code/.code={%
			\node[anchor=center] at (0.3cm,0cm) {#1};
		}
	},
}

\begin{tikzpicture}

\begin{axis}[%
width=4.44in,
height=2.533in,
at={(0.658in,0.5in)},
scale only axis,
xmin=0,
xmax=400,
xtick={0,50,100,150,200,250,300,350,400},
ymode=log,
ymin=0.0001,
ymax=3.57982399999988,
yminorticks=true,
axis background/.style={fill=white},
xmajorgrids,
ymajorgrids,
yminorgrids,
grid style={dotted},
ylabel = Relative Error,
xlabel = Basis Size $n$,
legend columns = 3,
legend style={
	font=\mystrut,
	legend cell align=left,
},
legend pos=north east,
every axis plot/.append style={thick},
]
\addlegendimage{legend image with text=POD-State}
\addlegendentry{}
\addlegendimage{legend image with text=POD-Indiv}
\addlegendentry{}
\addlegendimage{legend image with text=SVD-like}
\addlegendentry{}

\addplot [color=mycolor1, mark size=1.7pt, mark=triangle*, mark options={solid, fill=mycolor1, mycolor1}]
  table[row sep=crcr]{%
12	3.57982399999988\\
40	0.0724164099999967\\
80	0.0170019900000009\\
120	0.0114004199999995\\
160	0.0109368800000006\\
200	0.00879100499999946\\
240	0.00874029800000056\\
280	0.00699110799999971\\
320	0.00620347000000013\\
360	0.00590454599999995\\
400	0.00434961800000017\\
};
\addlegendentry{}

\addplot [color=mycolor2, mark size=2.5pt, mark=asterisk, mark options={solid, fill=mycolor2, mycolor2}]
  table[row sep=crcr]{%
12	0.734076700000039\\
40	0.307651100000019\\
80	0.0499447099999972\\
120	0.0222614399999999\\
160	0.0134775000000001\\
200	0.0116744699999993\\
240	0.0148310699999995\\
280	0.0102803800000004\\
320	0.0107592999999998\\
360	0.0088315419999995\\
400	0.00928286599999971\\
};
\addlegendentry{}

\addplot [color=mycolor3, mark size=2.5pt, mark=x, mark options={solid, fill=mycolor3, mycolor3}]
  table[row sep=crcr]{%
12	0.424234599999996\\
40	0.0588165900000019\\
80	0.00730455700000005\\
120	0.00282122599999998\\
160	0.00197608300000009\\
200	0.00201044099999993\\
240	0.00179845699999999\\
280	0.00181559699999988\\
320	0.00177896600000002\\
360	0.00167309999999999\\
400	0.00168081199999997\\
};
\addlegendentry{pH}


\addplot [color=mycolor1, dashed, mark size=1.7pt, mark=triangle*, mark options={solid, fill=mycolor1, mycolor1}]
  table[row sep=crcr]{%
12	1.01741400000001\\
40	0.0822529699999991\\
80	0.0196413399999995\\
120	0.0103936199999995\\
160	0.00829542099999969\\
200	0.00670506500000002\\
240	0.0061699690000004\\
280	0.0059911470000001\\
320	0.00537575600000008\\
360	0.00517551199999994\\
400	0.00345894899999994\\
};
\addlegendentry{}

\addplot [color=mycolor2, dashed, mark size=2.5pt, mark=asterisk, mark options={solid, fill=mycolor2, mycolor2}]
  table[row sep=crcr]{%
12	0.863229100000032\\
40	0.199991299999996\\
80	0.0341354100000013\\
120	0.0163867599999999\\
160	0.00744344699999998\\
200	0.00548381999999979\\
240	0.00417466999999984\\
280	0.00336119899999982\\
320	0.00310394900000006\\
360	0.00283566500000017\\
400	0.00270478899999984\\
};
\addlegendentry{}

\addplot [color=mycolor3, dashed, mark size=2.5pt, mark=x, mark options={solid, fill=mycolor3, mycolor3}]
  table[row sep=crcr]{%
12	0.406547099999974\\
40	0.0475498599999978\\
80	0.00736607499999987\\
120	0.00296762799999983\\
160	0.00185215800000009\\
200	0.00185347799999991\\
240	0.0016868330000001\\
280	0.00159102700000007\\
320	0.00147735800000009\\
360	0.0015831080000001\\
400	0.00159977600000003\\
};
\addlegendentry{Energy-Stable}

\addplot [color=mycolor1, dashdotted, mark size=1.7pt, mark=triangle*, mark options={solid, fill=mycolor1, mycolor1}]
  table[row sep=crcr]{%
12	0.309588100000006\\
40	0.0246833200000009\\
80	0.00654026600000037\\
120	0.00385561199999991\\
160	0.00305851099999989\\
200	0.00254792499999986\\
240	0.0021804050000001\\
280	0.00187807299999995\\
320	0.00153007900000001\\
360	0.00129895400000004\\
400	0.00107442200000004\\
};
\addlegendentry{}

\addplot [color=mycolor2, dashdotted, mark size=2.5pt, mark=asterisk, mark options={solid, fill=mycolor2, mycolor2}]
  table[row sep=crcr]{%
12	0.457640599999973\\
40	0.105289400000006\\
80	0.0138589399999993\\
120	0.0058272230000001\\
160	0.0029641769999999\\
200	0.00208616800000002\\
240	0.00172446800000003\\
280	0.00153594800000001\\
320	0.00141680999999999\\
360	0.00131887099999996\\
400	0.00123386600000001\\
};
\addlegendentry{}

\addplot [color=mycolor3, dashdotted, mark size=2.5pt, mark=x, mark options={solid, fill=mycolor3, mycolor3}]
  table[row sep=crcr]{%
12	0.353862200000004\\
40	0.0216502500000009\\
80	0.00368811700000002\\
120	0.00207661999999988\\
160	0.00164544499999991\\
200	0.00137762500000003\\
240	0.00113988599999997\\
280	0.000928123900000017\\
320	0.000763981800000025\\
360	0.000603919099999962\\
400	0.000462490100000012\\
};
\addlegendentry{Bestapprox}

\end{axis}
\end{tikzpicture}%
	\caption{relative error for three selected basis generation and projection methods in comparison with the best approximation}
	\label{fig:BestApproxVSMOR_MCMDS_selected}
\end{figure}
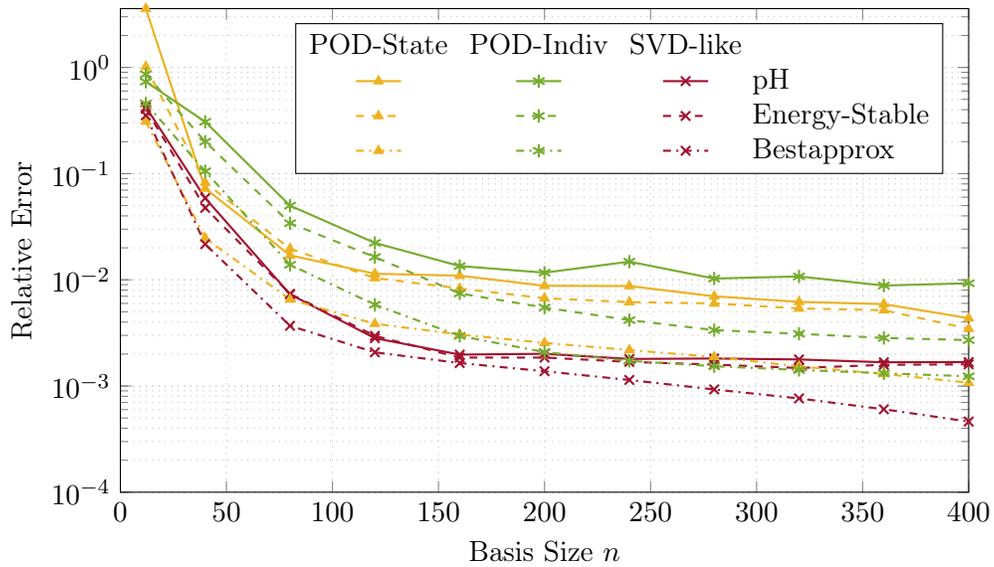

The three best-performing methods from the previous experiment, i.e. POD-State, SVD-like and POD-Individual, are compared more closely in Fig.~\ref{fig:BestApproxVSMOR_MCMDS_selected}. The plot shows the development of the error of the three methods over the basis size and examines on the one hand the pH and energy-stable projections, and on the other hand their behavior in relation to the best approximation with respect to the energy norm, i.e. the FOM solution is projected onto the subspace spanned by the corresponding basis matrix respecting the inner product pertaining to the energy norm \eqref{eq:BestApproximation}.
One can see that the results of pH and energy-stable behave in a similar way for POD-State and SVD-like but for POD-Individual, the energy-stable projection leads to much better error values. The deviation concerning the best approximation does not exceed one order of magnitude for all basis generation methods. The SVD-like method remains the closest to its potential best values for smaller reduced sizes with increasing distance for higher basis size. This is due to numerical reasons as the involved computation of the symplectic SVD-like decomposition is quite sensitive with respect to the condition of the snapshot matrix. A comparison of the basis generation methods uncovers that the SVD-like technique lies in the range of the best approximation error of the other methods and shows an improvement of approximately one order of magnitude for several reduction dimensions compared to the competitors. The overall behavior of all basis generation techniques (even the one not shown in the plot) shows that the errors decay rapidly up to a size of $n=120$ and stagnate starting from a size of approximately $n=200$.

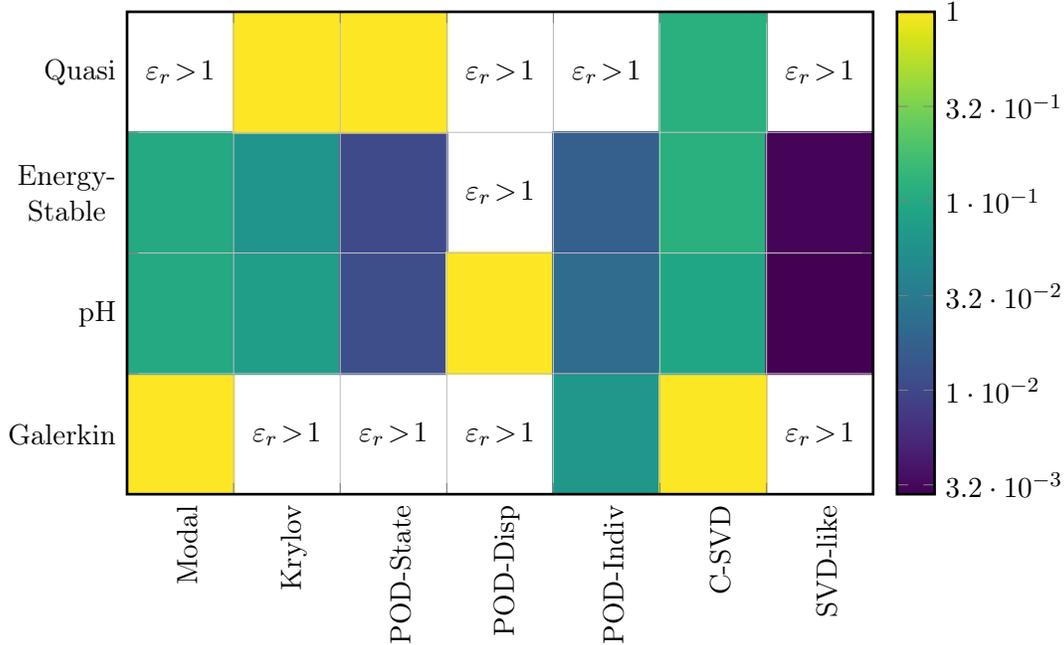
\begin{figure}[htp]
	\centering
	\pgfplotsset{compat=newest}


\pgfplotsset{every axis/.append style={line width = 1pt},%
    axis background/.style={fill=white},%
    width = 11.5cm,
	height = 8cm}

%



    
    \pgfmathdeclarefunction{lg10}{1}{%
    	\pgfmathparse{ln(#1)/ln(10)}%
    }
\begin{tikzpicture}
\newcommand\basisSize{120}
\newcommand\widthVar{12}
\newcommand\tresholdWhite{1}
\begin{axis}[
	xticklabels={Modal,Krylov,POD-State,POD-Disp,POD-Indiv,C-SVD,SVD-like},
	x tick label style ={rotate=90,yshift=-0.8cm},
	xtick={0.5,1.5,2.5,3.5,4.5,5.5,6.5},
	y tick label style ={yshift=0.8cm,align=center},
	yticklabels={Galerkin,pH,Energy-\\Stable,Quasi},
	ytick={0.5,1.5,2.5,3.5},
	axis on top,    
	colormap/viridis,
	colorbar,
	grid = major,
	point meta=explicit,
	mesh/ordering=y varies,
	mesh/cols=7,
	colorbar style={
	yticklabel={
	\pgfmathparse{10^\tick}\pgfmathprintnumber[precision=1]\pgfmathresult
}
	},
	xmin=0.5,
	xmax=7.5,
	ymin=0.5,
	ymax=4.5,
	]
	\addplot [matrix plot*] table [col sep=comma,point meta = lg10(\thisrow{errorValuesCSV})] {Figures/Results/SquareError_Impulse_Snap100_BasisSize\basisSize.csv};
	\node at (axis cs:2,1){$\varepsilon_{r}\!>\!\tresholdWhite$};
	\node at (axis cs:3,1){$\varepsilon_{r}\!>\!\tresholdWhite$};
	\node at (axis cs:4,1){$\varepsilon_{r}\!>\!\tresholdWhite$};
	\node at (axis cs:7,1){$\varepsilon_{r}\!>\!\tresholdWhite$};
	\node at (axis cs:1,4){$\varepsilon_{r}\!>\!\tresholdWhite$};
	\node at (axis cs:4,4){$\varepsilon_{r}\!>\!\tresholdWhite$};
	\node at (axis cs:5,4){$\varepsilon_{r}\!>\!\tresholdWhite$};
	\node at (axis cs:7,4){$\varepsilon_{r}\!>\!\tresholdWhite$};
	\node at (axis cs:4,3){$\varepsilon_{r}\!>\!\tresholdWhite$};
%
\end{axis}
\end{tikzpicture}
%
	\caption{relative errors for different basis and projection combinations for $n=120$}
	\label{fig:SquareErrorPlot}
\end{figure}

The heat map in Fig.~\ref{fig:SquareErrorPlot} shows the relative error for all combinations of basis generation techniques and projection methods that were conducted, cf.~\circledText{13}. The errors are shown for a basis size of $n=120$. If the relative error exceeds a value of $\varepsilon_r>1$ then the field entry is left blank. 

One can see that the Galerkin and Quasi-Galerkin projection perform the worst. These projections do not preserve the pH properties, especially the stability conditions, leading to instability and poor approximations. These projection methods are hence not suitable for the reduction of pH-systems. In general, the pH and energy-stable projection methods lead to the smallest error because of their adaption to the specific underlying pH structure of the system, which shows the importance of adapting the projection methods to the particular structure of the problem. 
The modal reduction, which is still used in various commercial MOR-packages, leads to high errors compared to the other combinations. The dynamics of the model cannot be described sufficiently accurate by only allowing to move in modal coordinates that belong to the lowest eigenfrequencies. The basis generated with the C-SVD and Krylov algorithms lead to moderate approximation errors which, in the case of the C-SVD, also cannot be improved by the transformation into the canonical structure of the matrix $\bm J$. The Krylov algorithm focuses on approximating the input-output behavior at certain frequencies, leading to higher error at different frequencies. The overall performance of the remaining data-based methods worked the best for the FSI dynamics of the guitar except for the POD-Displacement variant. The reason for the failure of the POD-Displacement has already been explained above. For the POD-State variant, the information contained in the snapshot matrix and extracted from this matrix via POD is sufficient for reproducing the coupled dynamics of the system. On the other hand, it can also be advantageous to consider the coordinates separately in the basis. The POD-Individual method outperforms the full-state alternative especially for the energy-stable projection.
This stems from the fact that the POD-Individual approach accounts for different scalings in the respective components. The SVD-like basis generation technique shows the best overall performance with error values as low as $\varepsilon_r = 3.2\cdot10^{-3}$. The explanation for this lies in the structure of the basis. The SVD-like approach is the only method that builds a non-orthogonal basis and can therefore adapt more flexibly to the system's dynamics. Special attention should be paid here to the fact that the SVD-like method was developed for purely Hamiltonian systems, whereas this study shows that it is also suitable for port-Hamiltonian systems.

It is worth mentioning that the information content in the snapshot matrices depends on the number of underlying trajectories which are randomly distributed in the parameter space. In the following study, trajectory numbers between five and 100 are investigated. The results of the associated relative errors are shown in Tab.~\ref{tab:ErrorOverSnapshots} and are extended by the methods based on system matrices for comparison.
It can be seen that the influence of a high number of trajectories is not as significant as for instance the size of the reduced system, see Fig.~\ref{fig:BestApproxVSMOR_MCMDS_selected}. All methods extract enough information for basis constructions even for small snapshot sizes. A snapshot matrix made of ten trajectories already shows a good trade-off between approximation error and matrix size, which drastically decreases the computational effort in the offline phase due to few evaluations of the full system.

\begin{table}
	\centering
	\caption{Relative error of various basis generation methods over the number of trajectories from the snapshot generation for the momentum formulation and $n=120$. One trajectory refers to an excitation with one specific frequency.}
	{\begin{tabular}{lrrrrrr} \toprule
			& \multicolumn{2}{l}{Number of Trajectories} \\ \cmidrule{2-7}
			Basis Generation & 5 & 10 & 25 & 50 & 75 & 100 \\ \midrule
			Modal & 0.0527 & \multicolumn{5}{c}{no data-based method, $\varepsilon_r$ for the sake of completeness} \\ 
			Krylov & 0.0293 & \multicolumn{5}{c}{no data-based method, $\varepsilon_r$ for the sake of completeness}  \\ 
			POD-State & 0.0051 & 0.0045 & 0.0047 & 0.0045 & 0.0043 & 0.0044 \\ 
			POD-Disp & 1 & 1 & 1 & 1 & 1 & 1 \\ 
			POD-Indiv &0.0118 &0.0104 &0.011&0.0114&0.0101& 0.0093 \\
			Complex SVD &0.0304 &0.0230 &0.0210&0.0218&0.0220&0.0207 \\
			SVD-like &0.0028 &0.0022&0.0021& 0.0018 & 0.0017&0.0017 \\ \bottomrule
	\end{tabular}}
	\label{tab:ErrorOverSnapshots}
\end{table}

The main goal of MOR is the gain of a speed-up compared to the simulation of the full-order model to make this model suitable for multi-query tasks or real-time scenarios. The speed-up values are determined for the pH projection method and averaged over all basis generation methods since the qualitative behaviour for the other methods was similar. The results are listed in Tab.~\ref{tab:SpeedUp}. The greatest speed-up of $540$ is obtained with a basis size of $n=12$. This choice is not recommendable since the error values are very high for this size. The speed-up values decrease down to $45$ for a basis size of $n=400$. Depending on the requirements for the ROM, one needs to decide which is the best trade-off between approximation quality and gained speed-up. Since the error decays fast up to a size of $n=120$ and still leads to speed-ups of more than 200, this reduction size was focused on in the previous discussion.

\begin{table}[htp]
	\centering
	\caption{Averaged speed-up values for different basis sizes for the pH projection}	
	{\begin{tabular}{lllllll}
			\toprule
			Basis Size & 12 & 40 & 80 & 160 & 320 & 400     \\ 
			Speed-Up & 540 & 455 & 340 & 180 & 68 & 45 		\\ \bottomrule
	\end{tabular}}
	\label{tab:SpeedUp}
\end{table}

\begin{figure}[htp]
	\centering
	\input{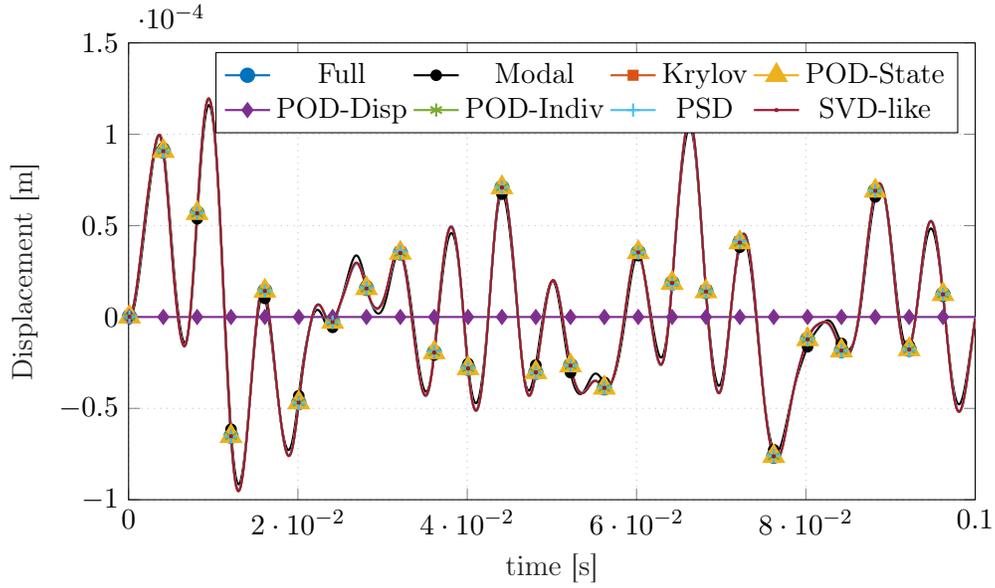}
	\caption{time simulation of the excitation node displacement for the full and reduced-order models for $n=120$ in the momentum formulation and pH projection}
	\label{fig:TimeSimulationRedTopPlate}
\end{figure}

\begin{figure}[htp]
	\centering
	\input{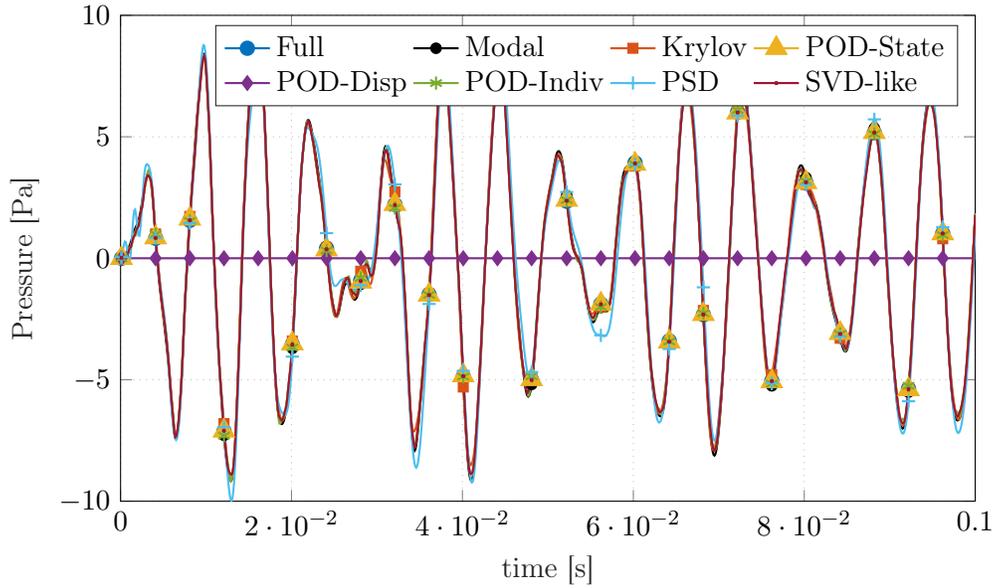}
	\caption{time simulation of the pressure integral over the sound hole for the full and reduced-order models for $n=120$ in the momentum formulation and pH projection}
	\label{fig:TimeSimulationRedSoundHole}
\end{figure}

In Fig.~\ref{fig:TimeSimulationRedTopPlate} and Fig.~\ref{fig:TimeSimulationRedSoundHole}, the displacement in the excitation node and the averaged pressure in the sound hole, cf.~Fig.~\ref{fig:FluidMesh}, are depicted. We can see in both cases, that all basis construction methods except for POD-Displacement capture the behaviour quite well. For the POD-Displacement, however, almost no energy seems to be introduced into the system resulting in no visible excitation in the entire FSI model.

\section{Conclusion and Outlook}
In this paper, we presented a port-Hamiltonian formulation of a fluid-structure interaction for the case of a classical guitar.
Many different model order reduction approaches combined with multiple different reduced basis construction methods as well as various FOM formulations were compared in order to study the effect of structure-preserving model reduction on the quality and behaviour of the reduced model. In particular, we were able to conclude that both structure-preserving model order reduction methods, i.e. the pH-preserving and the energy-stable methods clearly outperform their non-structure-preserving competitors. Furthermore, we can conclude that the reduced bases constructed via a symplectic SVD-like decomposition result in higher quality approximations for lower basis sizes.

Future emphasis will be placed on constructing suitable a-posteriori error estimators and extending the above model to allow for further parameter dependency in form of material parameters of the guitar body. \\


\acknowledgements{Supported by Deutsche Forschungsgemeinschaft (DFG, German Research Foundation)  Project Nos. 314733389, 455440338, and under Germany's Excellence Strategy - EXC 2075 – 390740016. We acknowledge the support by the Stuttgart Center for Simulation Science (SimTech).}

\bibliography{referenceGuitarPHMOR}
\bibliographystyle{abbrv}
\end{document}